\documentclass[11pt]{article}
\usepackage{a4wide}
\usepackage[dvips]{graphicx}
\usepackage{wasysym}
\usepackage{amsfonts}
\usepackage{amssymb}
\usepackage{epsfig}
\usepackage[all]{xy}
\newtheorem{defi}{\textbf{Definition}}[section]
\newtheorem{theo}[defi]{\textbf{Theorem}}
\newtheorem{lemma}[defi]{\textbf{Lemma}}
\newtheorem{prop}[defi]{\textbf{Proposition}}
\newtheorem{coro}[defi]{\textbf{Corollary}}

\newcommand{\Ang}{{\rm Ang}}
\newcommand{\MaxAng}{{\rm MaxAng}}
\newcommand{\Diff}{{\rm Diff}}
\newcommand{\Card}{{\rm Card}}
\title{Accidental Parabolics and Relatively Hyperbolic Groups}
\author{Fran\c{c}ois Dahmani \footnote{E-mail : dahmani@math.u-strasbg.fr}}
\date{}

\begin{document}
\maketitle

 {\footnotesize
 \begin{center} {\bf Abstract. } By constructing, in the relative case,
 objects analoguous to Rips and Sela's canonical representatives, we prove that
 the set of conjugacy classes of images by morphisms without accidental
 parabolic, of a finitely presented group in a relatively hyperbolic group, 
is finite.

\end{center}

}

\vskip .5cm

An important result of W.Thurston is :

\begin{theo}(\cite{Thur} 8.8.6)

 Let $S$ be any hyperbolic surface of finite area, and $N$ any geometrically finite
hyperbolic 3-manifold. There are only finitely many conjugacy classes of
subgroups $G\subset \pi_1(N)$ isomorphic to $\pi_1(S)$ by an isomorphism which
preserves parabolicity (in both directions).
\end{theo}

 It is attractive to try to formulate a group-theoretic analogue of this
statement : the problem is to find conditions such that the set of images of a
group $G$ in a group $\Gamma$ is finite up to conjugacy.

  If $\Gamma$ is word-hyperbolic and $G$ finitely presented, this has been the
object of works by M.Gromov (\cite{Gro} Theorem 5.3.C') and by T.Delzant
\cite{Del}, who proves the finiteness (up to conjugacy) of the set of images by
morphisms not factorizing through an amalgamation or an HNN extension over a
finite group.

 As a matter of fact, if a group $G$ splits as $A*_C B$ and maps to a group
$\Gamma$ such that the image of $C$ in $\Gamma$ has a large centralizer, then in
general, there are infinitely many conjugacy classes of images of $G$ in
$\Gamma$. Technically speaking, if $h$ is the considered map, one can conjugate
$h(A)$ by elements in the centralizer of $h(C)$, without modifying $h(B)$, hence
producing new conjugacy classes of images. A similar phenomenon happens with HNN
extensions.

 We are interested here in the images of a group in a relatively hyperbolic
 group (for example, a geometrically finite Kleinian group). Our result, Theorem
0.2, gives a condition similar to the one of Thurston, ruling out the bad
situation depicted above, and ensuring the expected finiteness.

 Relatively hyperbolic groups were introduced by M.Gromov in \cite{Gro}, and
studied by B.Farb \cite{Farb} and B.Bowditch \cite{Brel}, who gave different,
but equivalent, definitions  (see Definition \ref{def;relhyp} below, taken from
\cite{Brel}). In Farb's terminology, we are interested in ``relatively
hyperbolic groups with the property BCP''.  The main example is the class of
fundamental groups of geometrically finite  manifolds (or orbifolds) with
pinched negative curvature (see \cite{BGF}, see  also \cite{Farb} for the case
of finite volume manifolds). Sela's limit groups  are hyperbolic relative to
their maximal abelian non-cyclic subgroups, as shown in \cite{Dcomb}.

{\bf Definition :} We say that a morphism from a group in a relatively
 hyperbolic group $h : G\to \Gamma$ has an accidental parabolic \emph{either}
 if $h(G)$ is parabolic in $\Gamma$, \emph{or} if $h$ can be factorized
 through a non-trivial amalgamated free product $ \xymatrix{
        G \ar[r]^h  \ar@{->>}[rd]_f & \Gamma  \\
       & A*_C B \ar[u]
      }$
or HNN extension  $ \xymatrix{
        G \ar[r]^h  \ar@{->>}[rd]_f & \Gamma  \\
       & A*_C \ar[u] &
      }$
 where $f$ is surjective, and the image of $C$ is either finite or parabolic in $\Gamma$.

 We prove the theorem :

\begin{theo}

Let $G$ be a finitely presented group, and $\Gamma$ a relatively hyperbolic
 group. There are finitely many subgroups of $\Gamma$, up to conjugacy, that are images of
 $G$ in $\Gamma$ by a morphism without accidental parabolic.

\end{theo}

It would have been tempting to apply this to the mapping class group $Mod(S)$ of a
surface, which is known to be "relatively hyperbolic", after the study of
H.Masur and Y.Minsky of the complex of curves \cite{MM} (see also \cite{BCC}). 
If $B$ is the base
of a $S$-bundle, the study of homomorphisms $\pi_1(B) \to Mod(S)$ is important
because it is directly related to the geometric Shafarevich conjecture
(see the survey of C.McMullen \cite{McM}).
 Unfortunately, the relative hyperbolicity of the mapping class group is to be
understood in a weak sense : the property BCP, or equivalently the
\emph{fineness} (see Definition \ref{def;relhyp})
is not fulfilled.

Also note that Theorem 0.2 generalises Theorem 0.1 in the case of closed
surfaces: if a surface group $\pi_1(S)$ acts on a tree, an element associated
to a simple curve in $S$ fixes an edge. Therefore, if a morphism from
$\pi_1(S)$ to $\pi_1(N)$ (with notations of Theorem 0.1) has an accidental
parabolic, it sends a simple curve of the closed surface $S$ in a parabolic
subgroup of $\pi_1(N)$.

We will begin by introducing the definitions and the objects of the theory of
relatively hyperbolic groups. Then, in order to follow Delzant's idea in
\cite{Del}, we will generalize, in section 2, the construction of canonical
cylinders of Rips and Sela \cite{RS}  (Theorems \ref{theo;canonicalcylinders}
and  \ref{theo;slices}). The main difficulty comes from the fact that the
considered hyperbolic graph is no longer locally finite. Finally, we prove
Theorem 0.2 in section 3.

I would like to thank Thomas Delzant for the interesting discussions we had. I
am deeply grateful to Brian Bowditch, for his useful comments on this work.

\section{Relatively Hyperbolic Groups.}

\subsection{Definitions}

 A \emph{graph} is a set of \emph{vertices} with a set of \emph{edges}, which
are pairs of vertices. One can equip a graph with a metric where edges have
length 1. Thus this geometrical realization allows to consider geodesic,
quasi-geodesic and locally geodesic paths in a graph.
  A \emph{circuit} in a graph is a simple simplicial loop (without self intersection).
 In \cite{Brel}, B. Bowditch introduces fine graphs :

\begin{defi}(Fineness)\cite{Brel}

    A graph $\mathcal{K}$ is \emph{fine} if for all $L>0$, for all edge $e$, the
set of the circuits of length less than $L$, containing $e$ is finite. It is
\emph{uniformly fine} if this set has cardinality bounded above by a constant
depending only on $L$.
\end{defi}

We will use this definition  as a finiteness property of certain non-locally finite graphs.

\begin{defi}(Relatively Hyperbolic Groups)\label{def;relhyp}\cite{Brel}

   A group $\Gamma$ is \emph{hyperbolic relative to a family of subgroups}
$\mathcal{G}$, if it acts on a  hyperbolic and fine graph $\mathcal{K}$,  such
that  stabilizers of edges are finite,    the quotient $\Gamma \backslash
\mathcal{K}$ is a finite graph, and the stabilizers of the vertices of infinite
valence are exactly the elements of $\mathcal{G}$, and are
finitely generated.

    We will say that such a graph is \emph{associated} to the relatively
hyperbolic group $\Gamma$. A subgroup of an element of $\mathcal{G}$ is said to
be parabolic.

\end{defi}

  As there are finitely many orbits of edges, a graph associated to a relatively
hyperbolic group is uniformly fine.

 \subsection{Angles and Cones}

  As already explained in \cite{DY}, from which this section is partially
borrowed, angles and cones are useful tools for the study of fine graphs.

\begin{defi}(Angles)

Let $\mathcal{K}$ be a graph, and let $e_1=(v,v_1)$ and
$e_2=(v,v_2)$ be edges with one common vertex $v$. The angle
$\Ang_v(e_1,e_2)$, is the shortest length of the shortest of the paths from $v_1$
to $v_2$, in $\mathcal{K}\setminus \{v\}$. It is $+\infty$ if there is
no such path.

\end{defi}

The angle $\Ang_v(p,p')$ between two simple simplicial (oriented) paths $p$ and
$p'$ having a common vertex $v$ is the angle between their first edges after
this vertex.

If $p$ is a simple simplicial path, and $v$ one of its vertices,
$\Ang_v(p)$ is the angle between the consecutive edges of $p$ at
$v$, and its maximal angle $\MaxAng(p)$ is the maximal angle
between consecutive edges of $p$.

In the notation $\Ang_v(p,p')$, we will sometimes omit the
subscript if there is no ambiguity.

\begin{prop}(Three useful remarks)\label{prop;remarks}

 1.   When defined :
$\Ang_v(e_1,e_3)\leq \Ang_v(e_1,e_2)+\Ang_v(e_2,e_3)$.

 2.   If $\gamma$ is an isometry, $\Ang_v(e_1,e_2)= \Ang_{\gamma v}(\gamma.e_1,\gamma.e_2)$.

 3.  Any circuit (simple loop) of length $L \geq 2$ has a maximal angle less than $L-2$.

\end{prop}

 The first statement follows from the triangular inequality for the length
distance of $\mathcal{K}\setminus \{v\}$. The second statement is obvious.
 Finally, if $e_1=(v_1,v)$ and $e_2=(v,v_2)$ are two consecutive edges in the
circuit, the circuit itself gives a path of length $L-2$ from  $v_1$ and $v_2$
avoiding $v$. $\square$

\begin{defi}(Cones)

   Let $\mathcal{K}$ be a graph, $d>0$ and $\theta>0$.  Let $e$ be an edge, and
$v$ one of its vertices. The \emph{cone} centered at $(e,v)$, of radius $d$ and
angle $\theta$  is the set of vertices $w$ at distance less than $d$ from $v$
and such that there exists a geodesic $[v,w]$ satisfying the property that its
maximal angle and its angle with $e$ are less than $\theta$ :   $$
Cone_{d,\theta}(e,v)  =   \{w \, | \;  |w-v|\!\leq\! d,\,\MaxAng[v,
w]\!\leq\!\theta, \Ang_v(e,[v, w])\! \leq\! \theta \} $$

\end{defi}

\begin{prop}(Bounded angles imply local finiteness)

 Let $\mathcal{K}$ be a fine
 graph. Given an edge $e$ and a number $\theta>0$, there exists only finitely
many edges $e'$ adjacent to $e$ such that $\Ang(e,e')\leq \theta$. \end{prop}

{\it Proof : } There are only finitely many circuits shorter than $\theta$ containing $e$.
$\square$

\begin{coro}(Cones are finite)\label{cor;conesarefinite}

 In a fine graph, the cones are finite sets. If the graph is uniformly fine, the
cardinality of $Cone_{d,\theta}(e,v)$ can be bounded above by a function of $d$
and $\theta$. \end{coro}

 Consider a cone $Cone_{d,\theta}(e,v)$. We argue by induction on $d$. If $d=1$,
the result is given by the previous proposition. If $d>1$, we remark that
$Cone_{d,\theta}(e,v)$ is contained in the union of cones of angle $\theta$ and
radius $1$, centered at edges whose vertices are both in
$Cone_{(d-1),\theta}(e,v)$. If the latter is finite, the union is also finite.
$\square$

\begin{lemma}(Large angles in triangles)\label{lem;angles_triangles}

 Let $[x,y]$ and $[x,z]$ be geodesic segments in a $\delta$-hyperbolic graph,
and assume that $\Ang_x([x,y],[x,z])=\theta \geq 50\delta$. Then the
concatenation of the two segments is still a geodesic. Moreover $x$ belongs to
any geodesic segment $[y,z]$ and $\Ang_x([y,z]) \geq \theta -50\delta$.

\end{lemma}

{\it Proof : } Let $[y,z]$ be a geodesic, defining a triangle $(x,y,z)$, 
which is
$\delta$-thin. If both segments $[x,y]$ and $[x,z]$ are shorter 
than $10\delta$, then, $|y-z|\leq 20\delta$, and the total 
length of the edges of the triangle is less than $40\delta$,  The
third part of Proposition \ref{prop;remarks} proves that $x\in [y'',z'']$, 
and $Ang_x([y'',z'']) \geq \theta -50\delta $.

 Assume that $[x,y]$ is shorter than $10\delta$, and that 
$[x,z]$ is longer than $10\delta$. Let $z'$ be the point on $[x,z]$ at distance $13\delta$ from $x$. By triangular inequality, $|z'-y|\geq 3\delta$ and 
therefore, there exists a vertex $z''$ on $[y,z]$ at distance at most $\delta$ from $z'$. By triangular inequality, the segment $[z'',y]$ has length at most $24\delta$, and therefore, the loop $[x,z'][z',z''][z'',y][y,x]$ has length at most $48 \delta$. The segment $[z',z'']$  is at distance at least $12\delta$ from $x$, and therefore, it does not contain it. Again,  the
third part of Proposition \ref{prop;remarks} proves that $x\in [y'',z'']$, and
$Ang_x([y'',z'']) \geq \theta -50\delta $.

 Assume now that both  $[x,y]$ and $[x,z]$ are longer than $10\delta$.  
We consider the vertices $y'$ and $z'$ on $[x,y]$ and $[x,z]$ 
located at distance $10 \delta$ from $x$. 
If there is a path of length less than $3\delta$ between $y'$ and $z'$, it
cannot contain $x$, and therefore, it would contradict that the
angle at $x$ is greater than $50 \delta$. Therefore $y'$ and $z'$  are not
$3\delta$-close to each other, thus, they are $\delta$-close to the segment
$[y,z]$, and we set $y''$, respectively $z''$, in $[y,z]$ at distance less than
$\delta$ from $x'$, respectively $y'$. Consider the loop
$[x,y'][y',y''][y'',z''][z'',z'] [z',x]$. Its length is  less than $(2\times
10\delta+2\delta)\times 2  \leq 50\delta$, and it contains $x$. The segments 
$[y',y'']$ and $[z',z'']$ are at distance at least $9\delta$ from $x$, so that they do not contain $x$. Here again, the
third part of Proposition \ref{prop;remarks} proves that $x\in [y'',z'']$, and
$Ang_x([y'',z'']) \geq \theta -50\delta $. $\square$

\begin{lemma}\label{lem;circuitcone}(Cones and circuits)

  Let $e$ be an edge of a graph,
  and $w$ a vertex that lies in a circuit containing $e$ and of length less than
$L$. Then $w \in Cone_{L,L}(e,v)$. \end{lemma}

{\it Proof : } Let $C$ be the considered circuit, and let $g$ be a geodesic segment between
$v$ and $w$. The concatenation of $g$ and one of the two paths in $C$ from $w$
to $v$ is a loop. Hence, one has two loops containing $g$, one of them
containing $e$, one not, and both of length less than $L$. If $g$ has an angle
greater than $L$, then the corresponding vertex would not be in a sub-circuit of
each of the two loops, and therefore, the circuit $C$ would pass through this
point twice, which contradicts the definition of circuit. For the same reason
the angle between $e$ and $g$ is less than $L$, and therefore,   $w \in
Cone_{L,L}(e,v)$. $\square$

\begin{defi}
 Let $\Lambda$ be a number. A $\Lambda$-quasi-geodesic in a metric space $X$ is
a path $q:[a,b] \to X$ such that for all $x$ and $y$, $\frac{|x-y|}{\Lambda}
\leq {\rm dist}(q(x),q(y)) \leq \Lambda |x-y|$. \end{defi}

\begin{prop}(Conical stability of quasi-geodesics)\label{prop;finesseconique}

Let $\Lambda$ be a positive number. 
In a $\delta$-hyperbolic graph $\mathcal{K}$,
  let $g:[a,b] \to \mathcal{K}$ be a  geodesic segment, and let $q:[a,b] \to
\mathcal{K}$ be a $\Lambda$-quasi-geodesic with $|q(a)- g(a)| \leq r$ and
$|q(b)- g(b)| \leq r$, for $r\leq 10\delta$. Let $w$ be a vertex in $q$ at
distance at least $2r$ from the ends. Then there exists a constant
$N_{\Lambda,\delta}$ depending only on $\Lambda$, and $\delta$, and there exists
an edge $e$ in $g$, such that $w \in
Cone_{N_{\Lambda,\delta},N_{\Lambda,\delta}}(e,v)$.

\end{prop}

{\it Proof : } It is a classical fact  
(\cite{Gro}, 7.2 A, \cite{CDP}, \cite{GH}) that there 
exists a number $D(\Lambda, \delta)$ such that $q$ 
remains at a distance less than $D(\Lambda,\delta)$ from the segment, for a
certain constant $D(\Lambda,\delta)$.  We consider the loop starting at $w$,
consisting of five part : a subsegment $[w,w_1]$ of $q$, of length less than
$10.D(\Lambda,\delta)$, and strictly less if and only if $w_1=q(b)$, a segment
$[w_1,w_2]$  of length less than $D(\Lambda,\delta)$ and where $w_2 \in g$ (we
call it a transition),  a subsegment $[w_2,w_3]$ of $g$ of length  less than
$20.D(\Lambda,\delta)$ (strictly less if and only if  $w_3 = g(a)$), then again
a transition from $w_3$ to $q$  shorter than  $D(\Lambda,\delta)$, and then a
subsegment of $q$ to $w$.
As, in any case
$w$ is sufficiently far from the transitions, with respect to their length, it
does not belong to them, and this loop contains a sub-circuit shorter than
$25\Lambda D(\Lambda,\delta)$, containing $w$ and an edge of $g$.  Lemma
\ref{lem;circuitcone} gives the result. $\square$

\begin{figure}
\begin{center}
\input{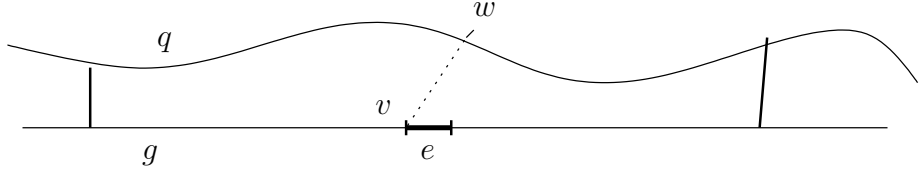}
 \caption{Quasi-geodesics stay in cones centered on the geodesic, Proposition
\ref{prop;finesseconique}}
\end{center}
\end{figure}

\section{Canonical cylinders for a family of triangles}

 In the following, $\mathcal{K}$ is a graph associated to a relatively
hyperbolic group $\Gamma$, and is $\delta$-hyperbolic. We choose a base point
$p$ in $\mathcal{K}$. We assume, without loss of generality, that $\delta$ is an integer greater than $1$.

The aim of this section is, given a finite family $F$
 of elements of $\Gamma$, to find a finite set  (a \emph{cylinder})  
around each
segment $[p, \gamma p]$ with $\gamma \in F\cup F^{-1}$. This construction will
be such that for all $\alpha, \beta, \gamma$ in $F\cup F^{-1}$ that satisfy the
equation $(\alpha \beta \gamma =1)$,  the three cylinders around $[p, \alpha
p]$, $[\alpha p, \alpha \beta p]=\alpha [p, \beta p]$ and $[p,
\gamma^{-1}p]=[\alpha \beta \gamma p, \alpha \beta p]$, coincide pairwise on
large subsets around the vertices $p$, $\alpha p$ and $\alpha \beta p$ (see
Theorem \ref{theo;canonicalcylinders}).

 Our approach  is similar to the original one in \cite{RS}. However, let us
emphasize that Rips and Sela use the fact that the balls in Cayley graphs are
finite. In the graph we are working on, it is not the case.

\subsection{Coarse piecewise geodesics}

We choose some constants :  $\lambda = 1000 \delta$, $\epsilon=N_{\lambda, \delta}$ and 
  $\mu =(100 \epsilon +  \lambda^2)\times 40 \lambda$ 
 ($N_{\lambda, \delta}$ is as in Proposition \ref{prop;finesseconique}).
These constants will be useful for defining \emph{coarse piecewise geodesics}, 
in 
the sense of \cite{RS}. Roughly speaking, $\lambda$ will serve as constant for quasi-geodesics, $\mu$ will serve as constant for local geodesics, 
and $\epsilon$ will be 
the bound for lengths of \emph{bridges}.

A path $p$ is a $\mu$-local-geodesic if any subpath of $p$ of length $\mu$, is a geodesic.
A path $p$ is a  $L$-local-$\frac{\lambda}{2}$-quasi-geodesic  if
 any subpath of length at most $L$ 
 is a $\frac{\lambda}{2}$-quasi-geodesic.

\begin{defi}(Coarse piecewise geodesics)(\cite{RS} 2.1) \label{def;cpg}

 Let $l$ be an integer greater than $\mu$. A $l$-coarse-piecewise-geodesic in $\mathcal{K}$ is a
$40\lambda(\epsilon+100\lambda\delta)$-local-$\frac{\lambda}{2}$-quasi-geodesic $f:[a,b] \to \mathcal{K}$ together
with a subdivision of the segment $[a,b]$, 
$a=c_1\leq d_1 \leq c_2 \dots \leq d_n =b$ such that 

\begin{itemize}

\item $f([c_i, d_i])$ is a $\mu$-local geodesic,   
\item $ \forall i,
2 \leq  i \leq (n-1), \;  length(f([c_i, d_i])) \geq l$ and $\forall i, \,
length(f[d_i, c_{i+1}])\leq \epsilon $,
\item $f([a,b])$ is included in the $2\epsilon$-neighborhood of a geodesic segment $[f(a),f(b)]$.
\end{itemize}

  In this case, we say that  $f|_{[c_i, d_i]}$ is a sub-local-geodesic, and
$f|_{[d_i, c_{i+1}]}$ is a bridge.

\end{defi}

\emph{Remark 1 : }  Any $l$-coarse-piecewise-geodesic is a $\lambda$-quasi-godesic. This follows from \cite{Gro} 7.2B, where it is stated that any 
$1000\delta \frac{\lambda}{2}$-local-$\frac{\lambda}{2}$-quasi-geodesic is a $\lambda$-quasi-geodesic. We also give, in appendix of this paper, a simple proof using the third point of the definition.

\emph{Remark 2 : } If $f:[a,b]\to \mathcal{K}$ is a coarse-piecewise-geodesic,
then for all $a'$ and $b'$ such that  $a\leq a'< b' \leq b$, the path
$f_{[a',b']}$ is a coarse-piecewise-geodesic. Indeed the induced subdivision satisfies the two first points of the definition (note that there is no length condition for the first and the last sub-local geodesic), and the third point is satisfied since $f$ (and therefore $f|_{[a',b']}$) is a $\lambda$-quasi-geodesic, and Proposition \ref{prop;finesseconique} applies.

\begin{lemma}(Re-routing coarse piecewise geodesics)
\label{lem;rerouting}

Let $l\geq \mu$ be a number, and $f:[a,b]\to \mathcal{K}$ be a 
$l$-coarse-piecewise-geodesic. 
Consider  a sub-local-geodesic $f|_{[c, d]}$, and 
$s \in f([c,d])$,  with  the additional requirement that the subpath of  
$f({[c, d]})$
from $f(c)$ to $s$ has length more than $l+2\epsilon $.
 Let $g$ be a geodesic segment between $f(a)$ and $f(b)$.
 Let $s''$ be a closest point to
$s$ on $g$. Let $s'$ be a closest point to $s''$ on $f_{[c, d]}$. 

We choose
$[s',s'']$ a geodesic segment between $s'$ and $s''$. We note $[s'',f(b)]$ a subsegment of $g$ between $s''$ and $f(b)$. We note $(f(a), s')$ the image by $f$ of the real segment $[a, f^{-1}(s')]$.

 Then the
path $\tilde{f}= (f(a), s') [s',s''] [s'',f(b)]$ is 
a $l$-coarse-piecewise-geodesic. 
We say that $f$ is re-routed into $\tilde{f}$.

\end{lemma}

{\it Proof : } 
We define the parameterisation of $\tilde{f}$ on the real 
segment $[a,b']$, to coincide with the one of $f$ on $[a,f^{-1} (s')]$, 
and to be the arc length parametrisation on $[f^{-1} (s'), b']$.  
Let $a=c_1< d_1\leq c_2<  \dots \leq c_m<d_m=b$ be the subdivision of $[a,b]$ associated to $f$, and let us say that $c=c_n$.  We define the 
subdivision of $[a,b']$, $a=c'_1< d_1\leq c'_2\leq  \dots \leq c'_n<d'_n\leq c'_{n+1} < d'_{n+1}=b'$ as coinciding with the one of the 
coarse-piecewise-geodesic $f$ until $c_n=c$, 
and such that $d_n = s'$, $c_{n+1} =s''$, and $d_{n+1}=b'$. 
It is immediate from similar property for $f$, 
that any restriction $f|_{[c_i,d_i]}$ is a 
$\mu$-local geodesic, that for all $i\in [2,n-1]$, 
$length(f([c_i,d_i]))\geq l$. and that for all $i\in [1,n-1]$, 
$length(f([d_i,c_{i+1}]))\leq \epsilon$. 
 We know that $f$ is a $\lambda$-quasi-geodesic, therefore by \ref{prop;finesseconique}, it stays $\epsilon$-close to $g$, hence, by triangular inequality, $\tilde{f}$ also stays $\epsilon$-close to $g$.

We have to show that  $length(f([d_n,c_{n+1}]))\leq \epsilon$, that $length(f([c_n,d_n]))\geq l$ and that $\tilde{f}$ is a $40\lambda(\epsilon + 100 \lambda \delta)$-local-$\frac{\lambda}{2}$-quasi-geodesic.

As $|s-s''|\leq \epsilon$, the segment  $[s',s'']$ has length less than $\epsilon$, which was the first requirement.

Therefore, $|s-s'| \leq 2\epsilon$, and, as we assumed that  the length of $f$
from  $f(c)$ to  $s$ is greater than $l+2\epsilon $
, the sub-local-geodesic 
of $\tilde{f}$ between $f(c)$ and $s'$ is longer than
$l$ in this case, which was the second requirement. 

We need to
prove that $\tilde{f}$ is a $40\lambda(\epsilon + 100 \lambda \delta)$-local-$\frac{\lambda}{2}$-quasi-geodesic. In
other words, we have to  show that any of its subpath of length less than 
 $40\lambda(\epsilon + 100 \lambda \delta)$
is a $\frac{\lambda}{2}$-quasi-geodesic.  Let $p$ be such a
subpath. If it is contained in the subpath of $\tilde{f}$ coinciding with
$f([a,b])$, by assumption on $f$ it is a
$\frac{\lambda}{2}$-quasi-geodesic. If it is contained in the subpath of
$\tilde{f}$ coinciding with $g$ it is a geodesic segment. If $p$ fails to
satisfy both conditions above, then it intersects $\rho$, and therefore is contained in a subpath of length at most $40\lambda(\epsilon + 100 \lambda \delta)- +2\epsilon$ that contains $\rho$. We give some
notations : let $x$ and $y$ be the ends of this subpath. 
As $\mu \geq 40\lambda(\epsilon + 100 \lambda \delta)
  +4\epsilon$, the subsegment $[x,s']$ of $f$ is a
geodesic. 
 The segment $[s'',f(b)]$ is included in
$g$ and therefore it is a geodesic segment, and it contains $y$. If  the length
of $p$ is less than $\frac{\lambda}{2}=500\delta$, there is nothing to prove.

It
is now enough to prove that for all such subpath $p$ containing $\rho$, 
of length more than $500\delta$, the distance $|x-y|$ between 
the ends $x$ and $y$ of $p$,
is superior to $\frac{1}{\lambda}\times (|x-s'|+|s'-s''|+|s''-y|)$.

 As the point $s'$ is the closest point to $s''$ in  $[x,s]$, by hyperbolicity,
we have $|x-s'|+|s'-s''| \leq |x-s''| +5\delta$.

Consider a point $u$ of the sub-local-geodesic $f([c,d])$ that is between $f(c)$
and $x$, and at distance $\mu/2$ from $x$. As  $l \ge \mu \ge 40\lambda(\epsilon + 100 \lambda \delta)+2\epsilon \ge |x-s'| $, it is possible to find such a point. Note
that the subpath $[u,s']$ of $f$ is of length at most $\mu$ and therefore
is a geodesic segment. 
Moreover, by Proposition \ref{prop;finesseconique}, there
is a point $v$ on $g$ such that  $|u-v|\leq \epsilon$. As the Gromov
product $(v\cdot y)_{s"}$ is equal to zero, and as $(v\cdot u)_{s"} \ge
\mu- 2\epsilon -10\delta \ge 100\delta$, by hyperbolicity, one has $(y\cdot
u)_{s"}\leq \delta$.

Similarily, $(u \cdot s")_x \leq 2\delta$, that is $(u\cdot x)_{s"} \geq |s"-x|
-5\delta$. There is the dichotomy : either $|s"-x| \le 20\delta$, hence
$|s"-y|\ge length(p) - 25\delta \geq |y-x| -25\delta$, and $(y\cdot
x)_{s"}\leq 45\delta$, or $|s"-x| \geq 20\delta$, and then $(u\cdot x)_{s"} \geq
 20\delta$, which together with $(y\cdot
u)_{s"}\leq \delta$, yelds $(y\cdot x)_{s"} \leq 2\delta$. In any case, one has $(y\cdot
x)_{s"}\leq 45\delta$. Then $|x-y| \geq |x-s"| + |s"-y| -45\delta$. We already
had  $|x-s'|+|s'-s''| \leq |x-s''| +5\delta$, which give : $|x-y|\geq
|x-s'|+|s'-s''|+|s''-y| - 50\delta$, and as $|x-s'|+|s'-s''|+|s''-y|$ was
assumed to be greater than $500\delta$, this gives the expected
$|x-y| \geq \frac{1}{\lambda}\times (|x-s'|+|s'-s''|+|s''-y|)$. This proves the
proposition. $\square$

\vskip .3cm
We will also need the following.

\begin{lemma}

 Let $[x,y]$ be a geodesic segment of $\mathcal{K}$, of length $L\geq
2\mu$. Let $s$ be on $[x,y]$ such that $|s-x|$ and $|s-y|$ are both
greater than $\frac{\mu}{2}$. Let $s'\in \mathcal{K}$ be at distance at most
$\delta$ from $s$ and $y'\in \mathcal{K}$ be at distance at most $\delta$ from
$y$. Let $s"$ be on $[x,y]$ such that $|s' - s"|$ is minimal. Then the path
$[x,s"] [s",s'] [s',y']$ is a $40\lambda(\epsilon + 100 \lambda \delta)$-local-$\frac{\lambda}{2}$-quasi-geodesic.


\end{lemma}

{\it Proof : } As in the previous lemma, it is enough to prove that for all subpath $p$ containing $[s",s']$, of length more than $500\delta = \frac{\lambda}{2}$, the distance $|x-y|$ between the ends $p_1$ and $p_2$ of $p$, is superior to
$\frac{1}{\lambda}\times (|p_1-s'|+|s'-s''|+|s''-p_2|)$.

 Let us assume that $|s''-p_2|\ge 25\delta$. By hyperbolicity, $p_2$ is $5\delta$-close to a point $w$ of $[s',y]$, and  $|s'-w| \geq |s''-p_2|- |p_2-w| -|s'-s"|$. Now $|p_1 -w| = |p_1 -s'| + |s' - w| \geq |p_1 -s'|+|s''-p_2|- |p_2-w| -|s'-s"|$. As $|p_1-p_2| \geq |p_1- w| - |w-p_2|$ we deduce that $|p_1-p_2|\geq|p_1 -s'|+|s''-p_2|- 2|p_2-w| -|s'-s"| \geq |p_1-s'|+|s''-p_2| +|s'-s"| -12\delta$, which is greater than $\frac{1}{1000\delta} \times(|p_1-s'|+|s'-s''|+|s''-p_2|)$, since $|p_1-s'|+|s'-s''|+|s''-p_2|$ is assumed to be greater than $500\delta$.

If $|s''-p_2|\leq 25\delta$, then $|p_1-p_2|\geq |p_1-s"| - |s''-p_2| \geq
|p_1-s'|+|s'-s''|+ |s''-p_2| -51\delta$, and the same conclusion holds.
$\square$

\begin{coro}(Re-routing to another point)\label{cor;reroutebis}

Let $l$ be a positive number. Let $f:[a,b]\to \mathcal{K}$ be a coarse-piecewise-geodesic whose last sub-local-geodesic $g$ is a geodesic segment of length at 
least $l+2\mu$. Let $[f(a),f(b)]$ be a geodesic segment of $\mathcal{K}$, and let $z$ be a point such that a geodesic segment $[f(a),z]$ passes at distance $\delta$ from $f(b)$. Then, there exists a $l$-coarse-piecewise-geodesic from $f(b)$ to $z$ coinciding with $f$ until the first point of $g$. 

\end{coro} 

{\it Proof : }  Let $f:[c,b] \to \mathcal{K}$ be an arc-length parametrization of the sub-local-geodesic $g$. One has $(b-c) \geq l+2\mu$. Let $x=f(c+l)$ and $y=f(b)$. 
Let $y'$ be a point of $[f(a),z]$ at distance less than $\delta$ from $y$. 



By the previous lemma, there exists points $s''$ on $g$, $s'$ on $[f(a),y']$
such that the path $f([a,c])[x,s"] [s",s'] [s',y']$ is a 
$40\lambda(\epsilon + 100 \lambda \delta)$-local-$\frac{\lambda}{2}$-quasi-geodesic satisfying the two first points of the definition of $l$-coarse-piecewise-geodesic. As $[s',y']$ is a subsegment of the geodesic segment $[s',z]$, the same is true for $f([a,c])[x,s"] [s",s'] [s',y']$.

 It remains to show that this paths stays $2\epsilon$-close to a geodesic segment $[f(a),z]$. Its first part from $f(a)$ to $s''$ is a subpath of $f$, hence it is  a $\lambda$-quasi-geodesic, therefore $\epsilon$-close to $[f(a),f(b)]$, and therefore, $(\epsilon + \delta)$-close to $[f(a),z]$. The second part $[s",s'] [s',z]$ is $\epsilon$-close to $[s',z]\subset [f(a),z]$ since $|s''-s'|\leq \epsilon$. This proves the claim. $\square$

     \subsection{Cylinders}

We now define the \emph{cylinders}, which are subsets of $\mathcal{K}$ associated 
to pairs of points.

\begin{defi}($l$-Cylinders)\cite{RS} \label{def;cylinders}

Let $l\in \Bbb{N}$. The $l$-cylinder of two points $x$ and  $y$ in
$\mathcal{K}$, denoted by $Cyl_l(x, y)$, 
is the set of the vertices $v$ lying on
a $l$-coarse-piecewise-geodesic from $x$ to $y$, with  the additional
requirement that $v$ is on a sub-local-geodesic $f|_{[c,d]}$ with distances
$|f(c) - v|\geq l$ if $f(c)\neq x$  and $|f(d)-v| \geq l$ if $f(d)\neq y$.
\end{defi}

Next lemma will assure that cylinders are finite sets, 
and stay close to geodesics.

\begin{lemma}(Cylinders are finite) \label{lem;cylinders_are_finite}

Given two points $x$ and $y$ in $\mathcal{K}$, and a constant $l$, any 
$l$-coarse-piecewise-geodesic from $x$ to $y$  remains in the union of the
cones of radius and angle $\epsilon$ centered in the edges of an arbitrary
geodesic segment $[x, y]$ 
(we call this union the $\epsilon$-conical-neighborhood of the segment). 

  The $l$-cylinder $Cyl_l(x,y)$  is contained in the union of the
cones of radius and angle $\epsilon$-conical-neighborhood of an arbitrary 
geodesic segment $[x, y]$.

The $l$-cylinder $Cyl_l(x,y)$  contains every geodesic between $x$ and $y$.
\end{lemma}

{\it Proof : } The second assertion is a consequence of the first one, itself being a consequence of  Proposition \ref{prop;finesseconique}  for $\Lambda = \lambda$, and $r=0$. To prove the third assertion it is sufficient to remark that every geodesic is a $l$-coarse-piecewise geodesic with only one sublocal geodesic, and no bridge. $\square$

\vskip .3cm

\vskip .3cm

\begin{lemma}(Equivariance)

 If a vertex $v$ is in $Cyl_l(x,y)$, then for all $\gamma$ in the group $\Gamma$, we have
$\gamma v \in Cyl_l(\gamma x, \gamma y)$
\end{lemma}

{\it Proof : } Multiplication on the left by $\gamma$ is an 
isometry of $\mathcal{K}$. $\square$

\subsection{Choosing a good constant $l$ for $l$-cylinders}

\begin{defi}(Channels)(\cite{RS} 4.1)

  Let $g=[v_1, v_2]$ be a geodesic segment in $\mathcal{K}$. A geodesic not
shorter than $|v_2- v_1|$ that stays in the union of the cones of radius and
angle $\epsilon$ centered in the edges of $g$ is a $(|v_2- v_1|)$-\emph{channel}
of $g$. \end{defi}

 As cones are finite (Corollary \ref{cor;conesarefinite}), the number of
different channels of a segment of length $L$ is bounded above by a constant
depending only on $L$. We note $Capa(L)$ such a bound.

 Let us recall the constants we fixed, and that are involved in the definition of 
coarse piecewise geodesics : $\mu =100 N_{\lambda,\delta}+  \lambda^2 $, with $\lambda =
1000\delta$.
For an integer $n$, we set
$\varphi(n) = 24(n+1)Capa(\mu)(2\epsilon+1)\epsilon $. For $1\leq i \leq
\frac{\varphi(n)}{2 \epsilon }$, we set now $l_i = 10\mu + 2i \epsilon $. Each  $l_i$ is
inferior to $\varphi(n)+10\mu$.

We denote by $B_r(x)$ the ball of $\mathcal{K}$ of center $x$ and radius $r$.

\begin{theo}\label{theo;canonicalcylinders}

 Let $F$ be a finite family of elements of $\Gamma$ ; we set $n = (2
\Card(F))^3$ where $\Card(F)$ is the cardinality of $F$. Let $p$ be a base point
in $\mathcal{K}$.

  There exists a number $l$ such that the $l$-cylinders satisfy : for all
$\alpha, \beta, \gamma $ in  $F\cup F^{-1}$ with $\alpha\beta\gamma =1$, in the
triangle $(x,y,z)=(p,\alpha p, \gamma^{-1} p)$ in $\mathcal{K}$, one has
$$ Cyl_l(x,y) \cap B_{R_{x,y,z}}(x) =  Cyl_l(x,z) \cap B_{R_{x,y,z}}(x)$$ (and
analogues permuting $x$, $y$ and $z$) where $R_{x,y,z} = (y\cdot z)_x -
4\times (11 \mu+\varphi(n))$, is the Gromov product in the triangle,
minus a constant.

\end{theo}

What is important in the theorem is not so much the value of $l$, but that the numbers
$ (y\cdot z)_x - R_{x,y,z}$ involved are bounded in terms of $n$ and of
$\mathcal{K}$ (namely, $\delta$ and  the cardinality of a
cone of radius and angle $\epsilon$). This bound does not depend on the family
$F$, although it does depend on its cardinality.

\begin{figure}
\begin{center}
\input{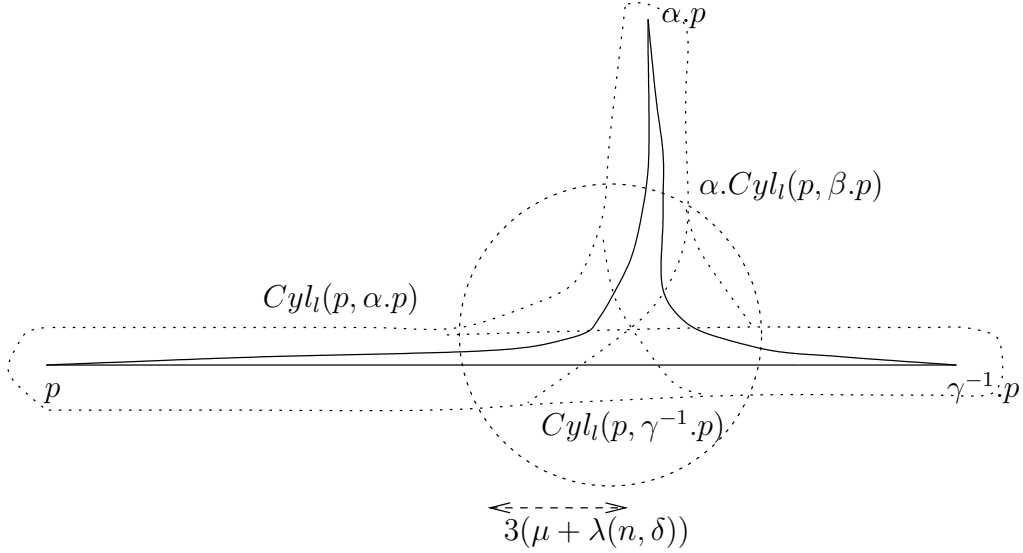}
\caption{Cylinders for a triangle, Theorem \ref{theo;canonicalcylinders}}
\end{center}
\end{figure}

{\it Proof : } We will find a correct constant $l$ among the $l_i$ previously defined. We have
$12 (n+1) Capa(\mu) (2\epsilon+1) $ different candidates. There are at most $n$
different triangles satisfying the condition, hence, we have a system 
of at most
$6n$ equations of the form  $Cyl_{l}(x,y)\cap
 B_{R_{x,y,z}}(x) \subset Cyl_{l}(x,z) \cap B_{R_{x,y,z}}(x)$. It is then enough to prove the next lemma.

\begin{lemma}  \label{lem;techniqueRS}

 Let $x,y,z$ be three points in $\mathcal{K}$. 
There are at most  $2Capa(\mu)(2\epsilon+1)$ different constants among the
$l_i$, $i=1..\frac{\varphi(n)}{2 \epsilon }$, such that   $Cyl_{l}(x,y)\cap
 B_{R_{x,y,z}}(x) \not{\subset} Cyl_{l}(x,z) \cap B_{R_{x,y,z}}(x)$.

\end{lemma}

{\it Proof : }  We argue by contradiction, assuming that
$(2Capa(\mu)(2\epsilon+1) +1)$ constants $l_i$ do not satisfy the equation  $Cyl_{l}(x,y)\cap
 B_{R_{x,y,z}}(x) \subset Cyl_{l}(x,z) \cap B_{R_{x,y,z}}(x)$. For
each of them, there is a vertex $v_i$ in one cylinder and not in the other :
there exists $\beta_i$, a $l_i$-coarse-piecewise-geodesic from  $x$ to $y$ containing $v_i$ as indicated in Definition \ref{def;cylinders}, but there is no such coarse-piecewise-geodesic 
from $x$ to $z$. 

By  Lemma \ref{lem;cylinders_are_finite}, each of the coarse-local-geodesics $\beta_i$ ($i=1..  \frac{\varphi(n)}{2 \epsilon }$) is contained in the $\epsilon$-conical-neighborhood (in the sense of Lemma \ref{lem;cylinders_are_finite}) of $[x,y]$.


  As $\epsilon \le \mu/2$, and $l_i \ge 10\mu$, 
each of the $\beta_i$ has a sub-local-geodesic passing through a $\mu$-channel
of a subsegment of $[x,y]$ starting at distance $ R_{x,y,z} + (\varphi
(n)+10\mu) $ from $x$ \emph{or} at distance 
$R_{x,y,z}  + (\varphi (n)+11\mu)$. 
There are less than $2 Capa(\mu)$ such channels.  
Therefore, there is a channel, denoted by
$Chan$, in which a sub-local geodesic
$\beta'_i=\beta_i|_{[c_i,d_i]}$ passes for at least $2\epsilon +2$ different indexes $i\in [1,\frac{\varphi(n)}{2 \epsilon}]$.  
Let us re-label $2\epsilon +2$ of these indexes : 
$i_1 <
i_2 < \dots < i_{2\epsilon +2}$.

For each $1 \leq j \leq 2\epsilon +2$, let $t_j \in [c_{i_j}, d_{i_j}]$ be the
instant where $\beta'_{i_j}(t_j)$ exits the channel $Chan$. Let us denote by
$r(\beta'_{i_j})$ the length of the path $\beta'_{i_j}([t_j,d_{i_j}])$, the
part of $\beta'_{i_j}$ after it leaves the channel $Chan$.
The discussion will hold on the respective possible values of the numbers
$r(\beta'_{i_j})$, for $1 \leq j \leq 2\epsilon +2$.

We now formulate and prove three claims.

{\bf Claim 1 : } For the any  $j\in [1,2\epsilon +2]$, one has  $r(\beta'_{i_j}) \leq l_{i_j}+2\epsilon$.
 
 Assume the contrary. Let $t^+_j>t_j$ be real number such that the length of  
$\beta'_{i_j}([t_j, t^+_j])$ equals  $l_{i_j}$.
 Then, by Lemma \ref{lem;rerouting},
$\beta_{i_j}$ can be rerouted into a $l_i$-coarse-piecewise-geodesic coinciding  with $\beta_{i_j}$ from $x$ to $\beta_{i_j} (t^+_j)$, coinciding with $[x,y]$ on a suffix starting at a point $3\epsilon$-close to $\beta_{i_j} (t^+_j)$, and ending at $y$. 
 Recall that $\beta_{i_j}(t_j)$ is the end point of the channel $Chan$. By triangular inequality, 
$|x- \beta_{i_j}(t_j)| \leq [R_{x,y,z} + (\varphi (n)+11\mu)+\epsilon +\mu]$. 
Therefore, $\beta_{i_j} (t^+_j)$ is at distance at most 
$[R_{x,y,z} + (\varphi (n)+11\mu)+\epsilon +\mu] + (\varphi
(n)+10\mu) + \epsilon \leq (y\cdot z)_x - 2\varphi (n) -22\mu$ from $x$.
Therefore, by Corollary \ref{cor;reroutebis}, at distance $l_{i_j}\leq (\varphi (n)+10\mu)$ after the bridge of the re-routing, it can be
rerouted into a $l_i$-coarse-piecewise-geodesic coinciding  with $\beta_{i_j}$ from $x$ to $\beta_{i_j} (t^+_j)$ that ends at  $z$. This
shows that $v_{i_1}$ is in the cylinders $Cyl_{i_1}(x,z)$, 
which contradicts our assumption, and proves the claim.

{\bf Claim 2 : } For any  two indices $i_j < i_k$,  one has $r(\beta'_{i_k}) <
r(\beta'_{i_j})$. 

If not, one could change $\beta_{i_j}$ just after $Chan$, 
into $\beta_{i_k}$ (it remains a $l_{i_j}$-coarse-piecewise-geodesic). 
On $\beta_{i_k}$, let  $\beta''_{i_k}$  be the 
sublocal geodesic following $\beta'_{i_k}$. It is on length 
$l_{i_k} \geq l_{i_j} +2\epsilon$. Let $t_k^+$ be real such that $\beta_{i_k}(t_k^+) $ is the point located on $\beta''_{i_k}$  at distance $l_{i_j} +2\epsilon$ from its beginning. 
 Recall (see paragraph above) that  
$|x- \beta_{i_j}(t_j)| \leq [R_{x,y,z} + (\varphi (n)+11\mu)+\epsilon +\mu]$. 
 By  Claim 1, we deduce that $\beta_{i_k}(t_k^+)$ is at distance at most 
$[R_{x,y,z} + (\varphi (n)+11\mu)+\epsilon +\mu] 
+ 2\times (\varphi (n)+10\mu) + 3 \epsilon 
\leq (y\cdot z)_x - 2\times (\varphi(n) + 10\mu + 2\epsilon)$ 
from $x$, and therefore, 
by Lemma \ref{lem;rerouting}, it is possible to reroute
the coarse-piecewise-geodesic, 
 into a $l_{i_j}$-coarse-piecewise-geodesic coinciding with $\beta_{i_j}$ until $Chan$, and coinciding with the suffix of $[x,y]$ starting at  $(y\cdot z)_x - 2\times (\varphi(n) + 10\mu$ from $x$, and  ending at $y$. 
Then, by Corollary \ref{cor;reroutebis}, it is possible to reroute it again on  into a $l_{i_j}$-coarse-piecewise-geodesic ending at $z$, which is a contradiction, as in Claim 1.


{\bf Claim 3 : } For all $i_k\leq 2\epsilon +2$, one has
$r(\beta'_{i_1})-r(\beta'_{i_k}) < 2\epsilon$.

 If not, we
could change $\beta_{i_k}$ just after $Chan$, by passing by $\beta'_{i_1}$, and
reroute it on $[x,y]$ before the end of  $\beta'_{i_1}$ 
(at distance $2\epsilon$
from the end). This again gives the same contradiction.

Now that the three claims have been proved, we can end the proof of Lemma 
\ref{lem;techniqueRS}.

 We see from the second claim that the  $2\epsilon +2$ numbers $r(\beta'_{i_j})$, for $j\leq
2\epsilon +2$, are all different, and, from the third claim, that they 
are all in an interval of $\mathbb{N}$ of
length $2\epsilon$ (hence containing $2\epsilon+1$ elements). This is a
contradiction. $\square$

\subsection{Decomposition of cylinders into slices}

 In this section we assume that the hypothesis of Theorem \ref{theo;canonicalcylinders} are fulfilled, and we choose $l$ a suitable constant as in the statement of this theorem. 
All considered cylinders will implicitly be $l$-cylinders.

Recall that there are only finitely many orbits of vertices of finite 
valence in $\mathcal{K}$, therefore, there exists a 
constant $\rho$ such that any pair of edges $(e,e')$ containing a vertex of finite valence $v$ makes an angle $\Ang_v(e,e')$ bounded above by $\rho$.

   Let $\Theta = \max \{10000 (D +\epsilon +\delta), \rho \}$, where $D$ is a constant such that
a $\lambda$-quasi-geodesic remains at distance $D$ from a geodesic in a
$\delta$-hyperbolic graph (here $\lambda = 1000 \delta$).

The decomposition into slices given by Rips and Sela in the hyperbolic case
(\cite{RS}) will not work properly here, because of large angles. Thus, we
choose a slightly different procedure.

\begin{defi}(Parabolic slices in a cylinder) \label{def;parab_slice}

 Let $a$ and $b$ two points in $\mathcal{K}$. In the cylinder $Cyl(a,b)$, a \emph{parabolic slice} is a singleton $\{v\}\subset
Cyl(a,b)$ such that there exists vertices $w$ and $w'$ in $Cyl(a,b)$, adjacent
to $v$ in $\mathcal{K}$ and such that  $\Ang_v((v,w),(v,w'))\geq \Theta$. The
angle of a parabolic slice is $Max_{w,w'\in Cyl}(\Ang_v ((v,w),(v,w'))$.
\end{defi}

Note that, since $\Theta \geq \rho$, any parabolic slice consists of a vertex of infinite valence. This justifies the name.

\begin{lemma}(Parabolic slice implies large angle on a geodesic segment) \label{lem;angularintrinsec}

Let $a$ and $b$ be two points in $\mathcal{K}$, and let $A$ be a number greater than $\Theta$. If $w$ and $w'$ are vertices in the cylinder $Cyl(a,b)$, such
that $|w-w'| \leq 50 \delta$, and if there exists $v$ on some geodesic $[w,w']$
such that $\Ang_v([v,w],[v,w'])=A$, then any geodesic segment
$[a,b]$ contains $v$, and $\Ang_v([a,b]) \geq A - 20D \geq A-\Theta$.

  If $\{v\}$ is a parabolic slice of a cylinder $Cyl(a,b)$, of angle $A$, then,
any geodesic segment $[a,b]$ contains $v$, and $\Ang_v([a,b]) \geq A - 20D
\geq A-\Theta$.

 \end{lemma}

{\it Proof : } The second assertion is an immediate corollary of the first one, and of Definition \ref{def;parab_slice}.

 Let $w$ and $w'$ be vertices in $Cyl(a,b)$, and $v$  
be such that  $|w-v| + |v-w'| =
|w-w'| \leq 50 \delta$ in $\mathcal{K}$, and such that $\Ang_v([v,w],[v,w'])=
A\geq \Theta$, for some geodesic segments $[v,w]$ and $[v,w']$.

As $w$ is in the cylinder $Cyl(a,b)$, there exists a $l$-coarse-piecewise-geodesic $f:[0, T] \to \mathcal{K}$, 
with $f(0)=a$ and $f(T)=b$, such that
 $f(s)=w$ for some $s\in [0,T]$, and such that 
$w$ is on a sub-local geodesic $f|_{[r,t]}$ of $f$, $|r-s|$ (resp.
$|s-t|$) being larger than $10\mu$, except if $r=0$ (resp ($t=T$).

As $f$ is a quasi-geodesic, at least one of the segments $f|_{[s,t]}$, and
$f|_{[r,s]}$ does not contain $v$. Let us assume that $f|_{[r,s]}$ does not
contain $v$.  
We set $s_1=\max\{0, s-3D\}$, and we choose  $x$ in a geodesic segment 
$[a,b]$ such
that the distance $|x-f(s_1)|$ is minimal (it is less than $D$, and
it is equal to $0$ if $s_1=0$). Let $[x, f(s_1)]$ be a geodesic segment.
 We claim that this segment does not contain $v$.  If $s_1=0$ the segment 
is exactly one
point equal to $a$, and it cannot be $v$ since $a$ is never a parabolic slice. 
If 
$s_1=s-3D$, the subpath $[f(s_1),w]$ of $f$ is included in a $\mu$-local geodesic, and is of length $3D<\mu$. Hence it is a geodesic, and therefore  $|f(s_1)-w| = 3D$. By triangular inequality, $|f(s_1)-v| \geq 3D -50\delta > |f(s_1)-x|$, and therefore, $[f(s_1),x]$ does not contain $v$, which is the claim.

Therefore there is a path $p$ from $w$ to
$x$ of length at most $4D$ not containing $v$.

We do the same construction for $w'$ : there exists $x'$ on $[a,b]$ and a path
$p'$ from $w'$ to $x'$ of length at most $4D$, not containing $v$. 
By triangular
inequality, $|x-x'| \leq 8D+ 50\delta \leq 9D$.

We now consider the path obtained by concatenation of $p$, $[x,x']$, and $p'$
(with reverse orientation). Its length is at most $17 D < A-50 \delta$. 
Therefore, the
segment $[x,x']$ must contain $v$, and the triangular inequality for angles
shows that $\Ang_v([x,x']) \geq A-17 D$. $\square$

\vskip .3cm

\begin{lemma}(Angles at the end of cylinders)\label{lem;pointes}

 Let $x\neq b$ be in $Cyl(a,b)$.
Then for all geodesic segments $[a,b]$ and $[x,b]$, $\Ang_b([x,b],[a,b]) \leq
14D$. \end{lemma}

{\it Proof : }  We distinguish two cases. First assume that $|x-b| \geq 3D$. 
We know that
there is a vertex $w$ on the segment $[a,b]$ such that $|w-x|
\leq D$. Therefore, in a geodesic triangle $(b,w,x)$, the segment $[b,x]$ and
$[b,w]$ remain $\delta$-close for a length at least $D \geq 10 \delta$.
Therefore, their angle at $b$ is less than $21\delta$, and it is less than
$14D$.

Secondly, assume that $|x-b| \le 3D$.
There is a coarse-piecewise-geodesic $f:[0,T] \to \mathcal{K}$ between $a$ and
$b$, containing $x$ on one of its sub-local geodesic. Let $t$ be such that
$f(t) =x$. Consider $t_1 = max\{0,t-3D\}$, and we choose $w\in[a,b]$ such
that the distance $|w-f(t_1)|$ is minimal (it is less than $D$ in any case, and
it is $0$ if $t_1=0$). Now we consider the path $p$ obtained by the
concatenation of a geodesic segment
$[w,f(t_1)]$ (of length at most $D$), of $f|_{[t_1,t]}$ (of length at most
$3D$), of a geodesic segment $[x,b]$ (of length at most $3D$), and of a
subsegment $[b,w]\subset [b,a]$ (of length at most $7D$ by triangular
inequality). As
$f$ is a quasi-geodesic, and $f(T)=b$, we deduce that $b$ is not on the path
$f|_{[t_1,t]}$. It is not on the segment $[w,f(t_1)]$ because $|w-f(t_1)| \leq
|f(t_1)-b|$. Therefore,  the path $p$ passes only once at the vertex $b$, and
therefore, $\Ang_b([x,b],[b,a])\leq 14D$.

We see that in any case, $\Ang_b([x,b],[b,a])\leq 14D$. $\square$

\vskip .3cm

\begin{lemma}(Angles in a cylinder)\label{lem;angle_in_cylinder}

 Let $[a,b]$ be a geodesic segment, such that for some vertex $v$ in $[a,b]$,
$\Ang_v([a,b]) > \Theta-20D$. Then, $Cyl(a,b) = Cyl(a,v)\cup Cyl(v,b)$.

In particular, if $\{v\}$ is a parabolic slice of $Cyl(a,b)$, then $Cyl(a,b) = Cyl(a,v)\cup Cyl(v,b)$.

Moreover, in such a case, $Cyl(a,v)\cap Cyl(v,b) = \{v\}$.

\end{lemma}

{\it Proof : }

Let $w$ be a point of $Cyl(a,b)$. There exists $f:[0,T]\to \mathcal{K}$  a $l$-coarse-piecewise-geodesic 
from $a$ to $b$ that contains $w=f(s)$ on one of its sub-local-geodesic, with the condition of Definition \ref{def;cylinders}. This coarse-piecewise-geodesic is a $\lambda$-quasi-geodesics, hence
stays $D$-close to the segment $[a,b]$. Hence, by an argument similar to Lemma
\ref{lem;angles_triangles}, any such coarse-piecewise-geodesic passes 
at the vertex $v$. Let $t$ be the real number such that $f(t)=v$. Then, by Remark 2, with the induced subdivision, $f|_{[0,t]}$ is a $l$-coarse-piecewise-geodesic from $a$ to $v$, and  $f|_{[t,T]}$ is a $l$-coarse-piecewise-geodesic from $v$ to $b$. Therefore, if $s\leq t$, we have that $w\in Cyl(a,v)$, and if $s\geq t$, then  $w\in Cyl(v,b)$.
This proves that $Cyl(a,b) \subset Cyl(a,v)\cup Cyl(v,b)$.

  Let us prove the other inclusion.  Let $w$ be a point of 
$Cyl(a,v)$. There exists  a
$l$-coarse-piecewise-geodesic $f :[0,T] \to \mathcal{K}$ from $a$ to $v$ containing $w$ on one of its sub-local-geodesic, with the condition of Definition \ref{def;cylinders}. 

Let $T' = T + |v-b|$, and let
$\tilde{f} : [0,T'] \to \mathcal{K}$ be as follows : $\tilde{f}|_{[0,T]} \equiv
f$, and $\tilde{f}(T+t)$ is the point of the given geodesic $[a,b]$ at distance
$T'-T-t$ from $b$. Let $f|_{[c,T]}$ be the last
sub-local geodesic of $f$, hence ending at $v$. 
Then $\tilde{f}|_{[c,T']}$ is still a
$\mu$-local-geodesic, by Lemma \ref{lem;angles_triangles}. Moreover, any subpath of length $1000\delta
\frac{\lambda}{2}  \leq \mu$ is a
$\lambda/2$-quasi-geodesic : either it is included in the path $f$, or in the
geodesic segment $[v,b]$, or it is the union of two geodesic segment that meet
at $v$ with an angle greater than $\Theta -20D$, and therefore is geodesic by Lemma \ref{lem;angles_triangles}. Finally, $\tilde{f}$ stays at distance $\epsilon$ from a geodesic segment $[a,b]$.  
Therefore, $\tilde{f}$ is a $l$-coarse-piecewise-geodesic from $a$ to
$b$, coincinding with $f$ between $a$ and $v$. This proves that the point $w$ is in $Cyl(a,b)$, and therefore, $Cyl(a,v) \subset Cyl(a,b)$.

 Similarly, by changing the role of $a$ and $b$, one has $Cyl(v,b) \subset Cyl(a,b)$ and therefore,  $Cyl(a,b) = Cyl(a,v)\cup Cyl(v,b)$.

The second assertion of the lemma is a consequence of  Lemma \ref{lem;angularintrinsec}.

Let us prove now that the intersection $Cyl(a,v)\cap Cyl(v,b)$ is $\{v\}$.
 Let $x$ be in the intersection $Cyl(a,v)\cap Cyl(v,b)$, and assume that $x\neq
v$. By Lemma \ref{lem;pointes}, $\Ang_v([x,v],[v,a])\leq 14D$.  Similarly, as
$x$ is also in $Cyl(v,b)$,  $\Ang_v([x,v],[v,b])\leq 14D$.
The triangular inequality for angles
(Proposition \ref{prop;remarks}) proves that $\Ang_v( [a,v], [v,b])$ is at most
$28 D$, and contradicts the assumption that it is greater than $\Theta-20D$. This prove
that $Cyl(a,v)\cap Cyl(v,b) = \{v\}$. $\square$

\vskip .3cm

The lemma we just proved allows us to consider unions of cylinders without
parabolic slice. This enables the contruction of regular slices, as in 
Rips and Sela  \cite{RS}.

 Let $Cyl(a,b)$ be a cylinder \emph{without parabolic slice}, and $x\in
Cyl(a,b)$. Following \cite{RS}, we define the set  $N^{(a,b)}_R(x)$ as follows : it is the set of all
the vertices $v\in Cyl(a,b)$ such that $|a-x|<|a-v|$, and such that
 $|x-v|>100\delta$.
 Here $R$ stands for ``right'', and $N^{(a,b)}_L(x)$
is similarly defined changing the condition $|a-x|<|a-v|$ into $|a-x|>|a-v|$.
As cylinders are finite, those sets are also finite.

\begin{defi}(Difference in cylinders without parabolic slice)\cite{RS}3.3

 Let $Cyl(a,b)$ be a cylinder without parabolic slice, and $x, y$ two points in
it. We define $\Diff_{a,b}(x,y) =\Card(N^{(a,b)}_L(x) \setminus N^{(a,b)}_L(y))
- \Card(N^{(a,b)}_L(y) \setminus N^{(a,b)}_L(x)) + \Card(N^{(a,b)}_R(y)
\setminus N^{(a,b)}_R(x)) - \Card(N^{(a,b)}_R(x) \setminus N^{(a,b)}_R(y))$,
where $\Card(X)$ is the cardinality of the set $X$. \end{defi}

Let us remark that this defines a cocycle (see \cite{RS}).

\begin{defi}(Regular slices in a cylinder without parabolic slice)

 Let $Cyl(a,b)$ be a cylinder without parabolic slice. An equivalence class in
$(Cyl(a,b) \setminus \{a,b\})$ for the equivalence relation
$(\Diff_{a,b}(x,y)=0)$ is called a regular slice of $Cyl(a,b)$.

\end{defi}

\vskip .3cm

 \emph{Ordering of slices.} We assign an index to each slice of $Cyl(a,b)$ as
follows.  Let $v_1,\dots, v_k$ be the
consecutive parabolic slices, ordered by their position on a geodesic segment
$[a,b]$. We set $S_0$ to be $\{a\}$. We define then $S_{j+1}$ to be the unique regular
slice of the cylinder $Cyl(a,v_1)$ such that $\Diff(S_{j},S_{j+1})$ is minimal.
If $S_j$ is the last slice in $Cyl(a,v_1)$, then the parabolic slice $\{v_1\}$
is labeled $S_{j+1}$. Then among the regular slices of a cylinder
$Cyl(v_i,v_{i+1})$, we define $S_{j+1}$ to be the (unique) slice such that
$\Diff(S_{j},S_{j+1})$ is minimal.
If $S_m$ is the
last regular slice of a cylinder $Cyl(v_i,v_{i+1})$ (for $i<k$), then the
parabolic slice $\{v_{i+1}\}$ is $S_{m+1}$. Finally we order the slices of the
last cylinder $Cyl(v_k,b)$ in the same way, and $\{b\}$ is the last slice (see
Figure 3).

\vskip .3cm

\begin{lemma} \label{lem;2delta}

 Let $a$ and $b$ be two points in  $\mathcal{K}$,  
and let $v$ be in $Cyl(a,b)$. Let $[a,b]$ be a geodesic segment. 
Then there exists $w\in [a,b]$ such that
$|w-v|\leq 2\delta$.
\end{lemma}

{\it Proof : }  The vertex $v$ is on a sub-local-geodesic of some
coarse-piecewise-geodesic $f$. By definition of the elements of cylinders, there
is a geodesic segment $[f(t_1),f(t_2)]$ containing $v$, such that, for $i=1,2$,
$f(t_i)$ is at distance at most $D$ of a
point $w_i \in [a,b]$, and such that \emph{either} $|v-f(t_i)| \geq 5D$ \emph{or} $f(t_i)$ equals to $a$ or $b$ (in this case, we choose $w_i$ to be $f(t_i)$).
The triangle $(w_1,f(t_1),f(t_2) )$ is  is $\delta$-thin. If the segment $[w_1,f(t_1)]$ is not reduced to a point, it remains at distance at least $4D$ from $v$, therefore $v$ is at distance at most $\delta$ from a point $v'$ of $[w_1,f(t_2)]$. Similarly, the triangle $(w_1,w_2,f(t_2) )$ is $\delta$-thin, and therefore, $v'$ is at distance at most $\delta$ from $[w_1,w_2]$.  Therefore, $v$ is at distance at most $2\delta$ from the segment $[w_1,w_2]$ included in $[a,b]$. $\square$

\vskip .3cm

\begin{lemma}\label{lem;noanglewithoutparabolic}

Let $Cyl(a,b)$ be a cylinder, and let $x$ and $y$ be two points in $Cyl(a,b)$.
Assume that there is a vertex $v$ in some geodesic segment $[x,y]$ such that
$\Ang_v([x,y])\geq 2\Theta$. Then, $\{v\}$ is a parabolic slice of $Cyl(a,b)$,
and if $x \in Cyl(a,v)$ then $y \in Cyl(v,b)$.

\end{lemma}

{\it Proof : }  If 
$|x-y|\leq 50 \delta$, the result is a consequence of Lemma \ref{lem;angularintrinsec}. 

If $|x-y|\geq 50
\delta$, let us parametrize the segment $[x,y]$ by arc length: $g:[0,L] \to
\mathcal{K}$, and let $g(t)=v$. By the  previous lemma, $x$ and $y$ are $2\delta$-close to a geodesic segment $[a,b]$. Let
$w$ and $w'$ be points in $[a,b]$ realizing this distance. 
 By hyperbolicity the triangles $(x,y,w)$ and $(x,w,w')$ are $\delta$-thin. Therefore, $g(t-5\delta)$ and $g(t+5\delta)$ is $2\delta$-close to the segment 
$[w,w']\subset [a,b]$.
Let us consider the path obtained by concatenation of geodesic segment :$[g(t+5\delta),y][y,w'][w',w][w,x][x,g(t-5\delta)]$, where the first and the last segments are included in $[x,y]$, and where $[w',w]\subset [b,a]$. By triangular inequality, its total length is at most $18\delta$. 
 We assumed that $\Ang_v([x,y])\geq 2\Theta$,  therefore, it must contain $v$, and the only possibility is that $v\in [w',w]\subset[b,a]$, moreover, $\Ang_v([w',w])\geq 2\Theta - 28 \delta$. 
By the third assertion of Lemma  \ref{lem;cylinders_are_finite}, 
every vertex of the segment $[a,b]$ is in the cylinder $Cyl([a,b])$. 
 Therefore, $\{v\}$ is a parabolic slice of $Cyl(a,b)$.

The second statement is a corollary of Lemma \ref{lem;pointes}. $\square$

\vskip .3cm

\begin{lemma}(Slices are small)\label{lem;slicesaresmall}

 Let $a$ and $b$ be two points of $\mathcal{K}$, and let $S$ be a slice of $Cyl(a,b)$. If $v$ and $v'$ are in $S$, then $|v-v'|\leq
200\delta$ and for all geodesic segment $[v,v']$, one has $\MaxAng([v,v'])\leq
2\Theta$.

\end{lemma}

{\it Proof : } If the slice $S$ is parabolic, then $v=v'$, and 
there is nothing to prove.
  We can  assume that the
slice $S$ is regular. We assume
without loss of generality that $|a-v|\leq|a-v'|$.

Suppose that $|v-v'|\geq 200\delta$. Let us prove that $N^{(a,b)}_L(v) \subset
N^{(a,b)}_L(v')$.

 By Lemma \ref{lem;2delta},
there is a vertex $w$ on a geodesic segment $[a,b]$ such that $|w-v| \leq
2\delta$, and similarily, there is
$w'$ on $[a,b]$ such that $|w'-v'| \leq 2\delta$. Note that, since $|v-v'| \geq
200\delta$, the distance $|w'-w|$ is at least $196 \delta$. The points $w$ and $w'$ are both on the segment $[a,b]$, therefore, if $|a-w |> |a-w'|$, then  $|a-w |\geq |a-w'|+196\delta$, and  $|a-v|>|a-v'|$, contradicting our assumption. Therefore, $|a-w |=|a-w|-|w-w'|\leq |a-w'| -196\delta$.

 Let $z$ be in $N^{(a,b)}_L(v)$. As it is an element of the cylinder $Cyl(a,b)$, by Lemma \ref{lem;2delta}, there
is an vertex $w_z$ of $[a,b]$ such that  $|z-w_z|\leq 2\delta$.  By definition, the vertex $z$
is at distance at least $100\delta$ from $v$, therefore $|w-w_z|\geq 96\delta$.

 Moreover, as $|a-z|\leq |a-v|$ and $|z-v| \geq 100\delta$, the vertex $w_z$ is
on the subsegment $[a, w]$ of $[a,b]$. 
Therefore, $|w_z-w'|=|w_z-w|+|w-w'|\geq 292\delta$.
This gives, by triangular inequality, $|z-v'|\geq 290 \delta$. Therefore,
$z$ is in $N^{(a,b)}_L(v')$.

Hence, we have $N^{(a,b)}_L(v) \subset
N^{(a,b)}_L(v')$ and similarly $N^{(a,b)}_R(v') \subset N^{(a,b)}_R(v)$.
 Moreover $N^{(a,b)}_L(v) \neq
N^{(a,b)}_L(v')$ (and similarly $N^{(a,b)}_R(v') \neq N^{(a,b)}_R(v)$),
 because $v'$ is in $N^{(a,b)}_L(v)$ and not in $N^{(a,b)}_L(v')$.
  Therefore, $\Diff_{a,b}(v,v') \neq 0$ which is a contradiction
since they both are in the same regular slice.

The bound on the maximal angle of a geodesic segment $[v,v']$ is a corollary of
the Lemma \ref{lem;noanglewithoutparabolic}: if  $\Ang_w([v,v'])\geq 2\Theta$
for some $w$, Lemma \ref{lem;noanglewithoutparabolic} implies that $v$ and $w$
are not in the same slice (not even in consecutive slices). $\square$

\vskip .3cm

\begin{figure}
\begin{center}
\input{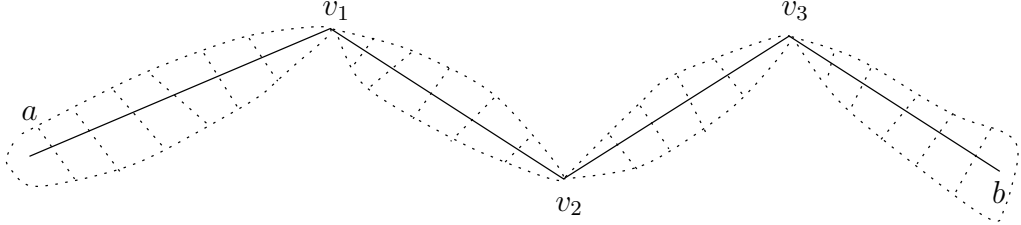}
\caption{Regular and parabolic slices in a cylinder}
\end{center}
\end{figure}

\begin{coro}(Consecutive slices are
close)\label{coro;slicesareclose}

Let $Cyl(a,b)$ be a cylinder, and let $S$ and $S'$ be two consecutive slices.
Let $v\in S$ and $v'\in S'$.

Then $|v-v'|\leq 1000 \delta$ and $\MaxAng([v,v']) \leq 2\Theta$.

\end{coro}

{\it  Proof : } The bound of the maximal angle is a consequence of Lemma
\ref{lem;noanglewithoutparabolic}: if there was such an angle there would be a
parabolic slice between $S$ and $S'$, hence, they would not be consecutive.

Assume that $|v-v'|\geq 1000 \delta$, and without loss of generality, $|a-v|\le
|a-v'|$ . By Lemma \ref{lem;2delta}, the points $v$ and $v'$ are $2\delta$-close to a geodesic segment $[a,b]$. Let $w$ be on
$[a,b]$, at distance at least $400\delta$ from $v$ and $v'$, and such that
$|a-v|\le |a-w| -200\delta \le |a-w| +200\delta \le |a-v'|$. By Lemma
\ref{lem;slicesaresmall}, $w$ is not in $S$ nor in $S'$, and as it is on a
geodesic segment $[a,b]$, it is in a slice. 
This slice is not before $S$ and not
after $S'$, therefore, $S$ and $S'$ are not consecutive. $\square$

\begin{lemma} (Locality of the regular slices)\label{lem;loc}

Let $Cyl(a,b)$, and $Cyl(a,c)$ be  cylinders without 
parabolic slices, and let $R$ be a number. 
Assume that  $Cyl(a,b) \cap B_{R}(a) = Cyl(a,c) \cap B_{R}(a)$, 
where $B_R(a)$ is the ball centered at $a$ of radius $R$. 
Then, any slice of $Cyl(a,b)$ included in
$B_{R-200\delta}(a)$ is a slice of $Cyl(a,c)$.

\end{lemma}

{\it Proof : } Let $S$ be a slice of $Cyl(a,b)$ and assume that $S$ is included in
$B_{R-200\delta}(a)$. Let $v$ be a point in $S$.
There exists $S'$, a slice of $Cyl(a,c)$ containing $v$. 
Let $v'$ be an arbitrary element of $S$. We claim that $v'$ is in $S'$. 

Let us compute
$\Diff_{a,c}(v,v')$. It is equal to   $\Card(N^{(a,c)}_L(v) \setminus
N^{(a,c)}_L(v')) - \Card(N^{(a,c)}_L(v') \setminus N^{(a,c)}_L(v)) +
\Card(N^{(a,c)}_R(v') \setminus N^{(a,c)}_R(v)) - \Card(N^{(a,c)}_R(v) \setminus
N^{(a,c)}_R(v'))$.

Note that $N^{(a,c)}_L(v) = N^{(a,b)}_L(v)$ and similarly for $v'$.
If $x$ is in $N^{(a,c)}_R(v') \setminus N^{(a,c)}_R(v)$, then it is
$100\delta$-close to $v$. Therefore, $x$ is in $Cyl(a,b)$, and it is  in
$N^{(a,b)}_R(v') \setminus N^{(a,b)}_R(v)$. Similarly the other inclusion
holds, and one has $N^{(a,c)}_R(v') \setminus N^{(a,c)}_R(v)=
N^{(a,b)}_R(v') \setminus N^{(a,b)}_R(v)$. 

Therefore,
$\Diff_{a,c}(v,v')=\Diff_{a,b}(v,v')$, and it is equal to $0$ since we assumed that  $v'\in S$. Therefore, $v' \in S'$, and  we deduce that $S' \subset S$.
Similarly, one has the other inclusion, and $S=S'$. This proves the lemma.
$\square$

\begin{theo}(Coincidence of the decomposition in slices)\label{theo;slices}

 With the notations of Theorem \ref{theo;canonicalcylinders}, let $(x,y,z) = (p,
\alpha p, \gamma^{-1}p)$  be a triangle in $\mathcal{K}$, such that $\alpha,
\beta, \gamma$ are in $F \cup F^{-1}$, and $\alpha \beta \gamma =1$.

The ordered slice decomposition of the cylinders is as follows.
$$
\begin{array}{lllllllllll}
 Cyl(x,y) & = &(S_1, & S_2,  & \dots, & S_k, &\, \mathcal{H}_z,\, & T_m,  & T_{m-1}, & \dots, & T_1)   \\
 Cyl(x,z) & = &(S_1, & S_2,  & \dots, & S_k, &\, \mathcal{H}_y,\, & V_p, & V_{p-1}, & \dots, & V_1)   \\
 Cyl(y,z) & = &(T_1, & T_2,  & \dots, & T_m, &\, \mathcal{H}_x,\, & V_p, & V_{p-1}, & \dots, & V_1),
\end{array}
$$

 where $S_1, \dots, S_k, T_1,\dots, T_m$ and  $V_1,\dots, V_p$ are slices and where each $\mathcal{H}_v$,
$(v=x,y,z)$ is a set of at most $10 \varphi(n)$ consecutive slices,
without parabolic slice of angle more than $3\Theta+ 100\delta$.

The sets  $\mathcal{H}_v$ are called the \emph{holes} of the slice decomposition.

\end{theo}

{\it Proof : } Consider the cylinders $Cyl(x,y)$ and $Cyl(x,z)$.
  By Theorem \ref{theo;canonicalcylinders}, they coincide in $B_{R_{x,y,z}}(x)$.
Therefore any parabolic slice of $Cyl(x,y)$ that is located in
$B_{R_{x,y,z}-2}(x)$ is also a parabolic slice of $Cyl(x,z)$, and similarly, permuting $x$ and $y$.

Let $\{v\}$ be the last common parabolic slice of these two cylinders, or $v=x$ if they have no common parabolic slice: $Cyl(x,y) = Cyl(x,v) \cup
Cyl(v,y)$ and $Cyl(x,z) = Cyl(x,v) \cup
Cyl(v,z)$, by Lemma \ref{lem;angle_in_cylinder}.

The ordered slices of the cylinders $Cyl(x,y)$ and $Cyl(x,z)$ obviously coincide
at least until the slice $\{v\}$.

Let $\{w\}$ be the first parabolic slice of $Cyl(x,y)$ after $\{v\}$, or $w=y$
if there is no such parabolic slice. Let $\{w'\}$ be the first parabolic slice
of $Cyl(x,y)$ after $\{v\}$, or $w'=z$ if there is none. By Theorem
\ref{theo;canonicalcylinders}, $Cyl(v,w) \cap B_{R_{x,y,z} -|x-v|} (v) =
Cyl(v,w') \cap B_{R_{x,y,z} -|x-v|} (v)$. These cylinders are without parabolic
slices. By Lemma \ref{lem;loc}, their regular slices lying in $B_{R_{x,y,z}
-|x-v| -200\delta} (v)$ coincide.

 In other words, the slice decomposition of  $Cyl(x,y)$ and  $Cyl(x,z)$ coincide
at least until their last common parabolic slice, and for all slices in
$B_{(R_{x,y,z}-200\delta)}(x)$. A similar statement holds for the other pairs of
cylinders.

It remains to prove that no hole contain a parabolic slice of angle 
greater than  $3\Theta+100\delta$. Let us consider such a parabolic slice $S=\{v\}$ in $Cyl(x,y)$. By Lemma
\ref{lem;angularintrinsec}, a segment $[x,y]$ has an angle greater than
$2\Theta+100\delta$ at the point $v$.  Therefore, one of the two segments
$[x,z]$ and $[z,y]$ (say $[x,z]$) has an angle greater than $\Theta$ at $v$, 
and we deduce that $\{v\}$ is a parabolic slice of $Cyl(x,z)$. 
As it is simultaneously a parabolic slice in  $Cyl(x,y)$ and in $Cyl(x,z)$, 
it is not in a hole. $\square$

       \section{Image of a group in a relatively hyperbolic group}

In this section we consider $\Gamma$ a relatively hyperbolic group with 
associated graph
$\mathcal{K}$, and $G$ a finitely presented group with a morphism $h : G \to
\Gamma$. We want to explain how to adapt Delzant's method, given for
hyperbolic group in \cite{Del}, to the relative case, in order to obtain an
analogue to Thurston's Theorem 0.1.

 For conveniance, we choose the graph  $\mathcal{K}$ with the four
following properties

It has a base point $p$ with trivial stabilizer. Its vertices are exactly the
infinite valence vertices and the elements of the orbit of $p$. It has no
pair of adjacent vertices of infinite valence. Finally,  for a
certain word metric on $\Gamma$, one has, for all $\gamma$ in $\Gamma$, for all
geodesic segment $[p,\gamma p]$ in $\mathcal{K}$, 

\begin{equation} \label{eq;majoration} 
|\gamma p - p|\times (\MaxAng([p,\gamma p]) +1) \geq |\gamma|.
\end{equation}

 It is possible to choose $\mathcal{K}$ satisfying these requirements : see for
example the coned-off graph of the Cayley graph in \cite{Farb}, 
where the angles
at the parabolic vertices are bounded by a word metric of the parabolic
subgroups, which are assumed to be finitely generated. To see that such the majoration (\ref{eq;majoration}) is fulfilled, it suffices to see that the distance beetween two points in the Cayley graph is bounded above by the length of a path of $\mathcal{K}$ between these two points plus the sum of the angles of this path at the vertices of infinite valence. The three other conditions are obvious.

\emph{Remark 3 : } In such a graph, a cylinder cannot have two consecutive
parabolic slices. Indeed, a geodesic segment between two parabolic slices
$\{v_1\}$ and $\{v_2\}$ must contain a vertex with trivial stabilizer, which
would belong to some regular slice of $Cyl(v_1,v_2)$.

\begin{defi}(Accidental parabolic)

 We say that the morphism $h : G \to \Gamma$ has an accidental  parabolic either
if  $h(G)$ is parabolic in $\Gamma$, or if there
 exists a non-trivial amalgamated free product $A*_C B$, or an HNN extension
$A*_C$,  and a factorization of $h$ : $ \xymatrix{         G \ar[r]^h
\ar@{->>}[rd]_f & \Gamma  \\        & A*_C B \ar[u]_{h'} &
      }$
 or  $ \xymatrix{
        G \ar[r]^h  \ar@{->>}[rd]_f & \Gamma  \\
       & A*_C  \ar[u]_{h'} &
      }$
 such that $f$ is surjective and the image of $C$ by $h'$ is a  finite, or
parabolic subgroup of $\Gamma$. \end{defi}

\begin{lemma}\label{lem;virtparab}
 If a subgroup $H$ of $\Gamma$ has a finite orbit in the graph $\mathcal{K}$,
then  either $H$ is finite or it is parabolic. \end{lemma}

{\it Proof : }  The subgroup $H$ has a subgroup of finite index $P$, 
fixing a point in
$\mathcal{K}$. Assume that $H$ is infinite, and not equal to $P$. 
As $P$ is also
infinite, it is parabolic, and the intersection of all its conjugates in $H$ is
infinite. But it is easily seen from fineness that the intersection of two
distinct conjugates of a maximal parabolic subgroup is finite in a relatively
hyperbolic group. Hence, $H$ is itself parabolic. $\square$

In the rest of this section, we prove the next theorem.

\begin{theo}\label{theo;acci}

  Let $G$ be a finitely presented group, and $\Gamma$ a relatively hyperbolic
group.  
There is a finite family of subgroups of $\Gamma$ such that the image of
$G$ by any morphism $h:G\to \Gamma$ without accidental parabolic is conjugated
to one of them.

\end{theo}

{\it Proof : } Let $h$ be a morphism $h:G\to \Gamma$. We will construct a factorisation of
$h$ through a certain graph of groups, and then we will deduce that either $h$
has an accidental parabolic, or $h(G)$ is conjugated to a subgroup of $\Gamma$
generated by small elements.

We choose a triangular presentation of $G$ : $G=<g_1,\dots,g_k | T_1, \dots,
 T_n>$ with $n$ relations which are words of three (or two) letters.
 This defines a Van Kampen polyhedron $P$ for $G$,  which consists of $n$
triangles and digons.

  Recall that the base point $p$ of the graph $\mathcal{K}$ associated to the
relatively hyperbolic group $\Gamma$, has trivial stabilizer.
 We consider the cylinders of the triangles, and their decomposition in slices
 obtained by the Theorems \ref{theo;canonicalcylinders} and \ref{theo;slices},
for the family
 $F=\{h(g_1), \dots, h(g_k)\} \subset \Gamma$ and the base point $p \in
 \mathcal{K}$.

 \subsection{The lamination $\Lambda$ on $P$.}

We now define a lamination on $P$, in two steps : first by choosing markings on the edges of $P$, and secondly by defining arcs in $P$ between these markings.

\subsubsection{Markings on the edges of $P$}

 For a generator $g_i$ of $G$, let $L^r_i$ be the number of regular slices of the
cylinder of $[p,h(g_i)p]$ in $\mathcal{K}$, and $L^p_i$, the number of its
 parabolic slices. Let $c_i$ the loop of the
  polyhedron $P$ canonically associated to $g_i$. Let $m_i^1, \dots m_i^{L^r_i +
2L^p_i}$ be $(L^r_i + 2L^p_i)$ points on $c_i$, such that, if $c_i (t):[0,1]
\to P$ is an arc-length parametrisation of $c_i$, one has $m_i^k = c_i( k
\frac{1}{L^r_i +
2L^p_i+1})$.
 We call them the \emph{markings} of the slice decomposition on
$c_i$.
To each marking of $c_i$ we associate a slice in the cylinder of
$[p,h(g_i)p]$ in $\mathcal{K}$ : $m_i^1$ is associated to the first slice ; if
$m_i^k$ is associated to a regular slice, $m_i^{k+1}$ is associated to the next
slice in the ordering, if
$m_i^k$ is associated to a parabolic slice, and if $m_i^{k-1}$ is associated to
another slice (or do not exist), then $m_i^{k+1}$ is associated to the same
slice than $m_i^k$ ; finally, if $m_i^k$ and $m_i^{k-1}$ are associated to the
same parabolic slice, then $m_i^{k+1}$ is associated to the next
slice in the ordering. Note that every regular slice has one marking on $c_i$
associated to it, and every parabolic slice has two markings.

\subsubsection{Regular arcs in a triangle (or a digon) of $P$}

The lamination $\Lambda$ is defined on $P$ by its intersection with each triangle or
digon $T$ in $P$.

 Consider a triangle $T$ (with an euclidean metric) of $P$, whose edges
$c_i,c_j,c_k$ correspond to the relation $g_{i}g_{j}g_k=1$ of  the presentation.

 Consider two markings $m_i^{r}$ of
$c_i$ and $m_j^{s}$ of  $c_j$, that are  associated to the same regular slice
in the cylinders of the triangle $(p,h(g_i) p, h(g^{-1}_k) p)$ in
 $\mathcal{K}$. The segment $[m_i^{r}, m_j^{s}]$ in $T$ is said to be a
\emph{regular arc}.

  Consider two consecutive markings, $m_i^{r}$ and $m_i^{r+1}$, of $c_i$,
associated to the same parabolic slice of $Cyl([p,h(g_i)p])$. There are three
possibilities. 

First, if the slice is not equal to a slice of any of the two other
cylinders (that is: if it is in a hole in the sense of Theorem \ref{theo;slices}), 
we do nothing. 

Secondly, if it is a slice of one, and only one, other
cylinder, say $Cyl([h(g_i)p, h(g_j)h(g_i)p])$, then there are two consecutive
markings $m_j^{s}$ and $m_j^{s+1}$ of $c_j$ associated to it. The segments
$[m_i^{r}, m_j^{s+1}]$ and $[m_i^{r+1}, m_j^{s}]$ are said to be also
\emph{regular arcs}.
Note that these two segments do not cross.

 Finally, if the
slice is a slice of $Cyl([p, h(g_j)p]$ and of $Cyl([h(g_j)p, h(g_i)p])$, there
are two consecutive markings $m_j^{s}$ and $m_j^{s+1}$ of $c_j$, and two
consecutive markings $m_k^{t}$ and $m_k^{t+1}$  of $c_k$, associated to it. The
three segments $[m_i^{r}, m_k^{t+1}]$, $[m_i^{r+1}, m_j^{s}]$, and $[m_j^{s+1},
m_k^{t}]$ are \emph{regular arcs}. These three segments do not cross each other.

We do similarly after cyclic permutations
of $i,j$ and $k$.
We denote by $\Lambda_r(T)$  the union of all the regular arcs in $T$.

\subsubsection{Singular arcs in a triangle (or a digon) of $P$}

 If the
slice decomposition of the triangle has a hole (in the sense of Theorem
\ref{theo;slices}), there are markings that are not in
regular arcs.  In such a case,  we add a
\emph{singular  point} $p_T$ in the component of $T\setminus \Lambda_r(T)$
containing these markings. For all marking $m$ not in  $\Lambda_r(T)$, the
segment $[m,p_T]$ is said to be a \emph{singular arc}. Let $\Lambda_s(T)$ be
the union of these singular arcs in $T$.

The lamination $\Lambda$ on $P$ is defined by : for all triangle or digon $T$
of $P$, $\Lambda \cap T = \Lambda_r(T)\cup \Lambda_s(T)$ (see figure 4).

\begin{figure}
\begin{center}
\input{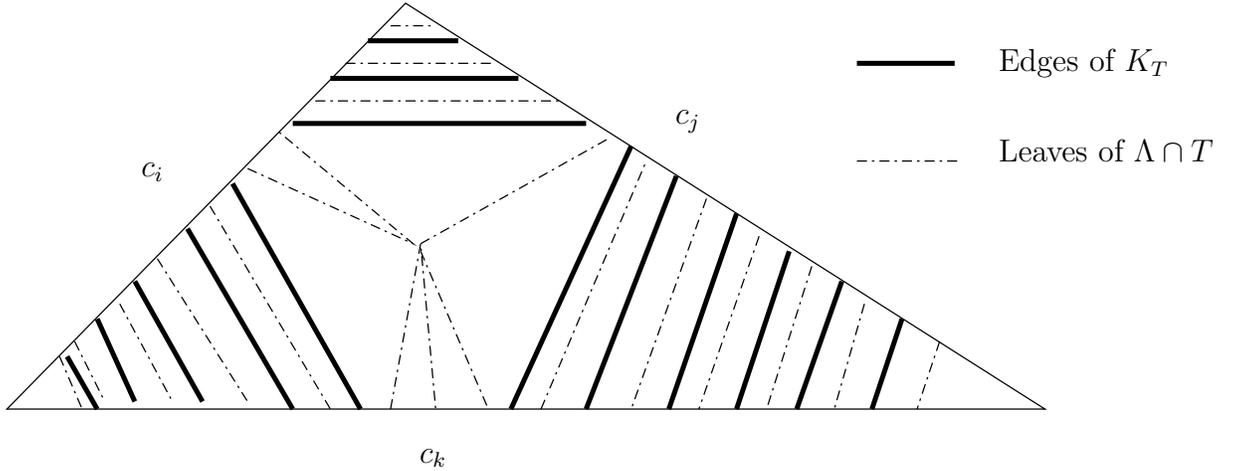}
\caption{The lamination $\Lambda \cap T$ and the graph $K_T$ in a triangle $T$
of $P$} \end{center}
\end{figure}

\subsection{Graph $K$ on $P$.}

In each triangle or digon $T$ of $P$, we draw a (disconnected) graph
$K_T$ satisfying :
each connected component of $T\setminus K_T$ contains one and only one leaf
of $\Lambda \cap T$, and its intersection with the edges $c_i$ of $T$, $K_T
\cap c_i$ consists of  the
vertices of $K_T$, moreover they are located on middles  of consecutive markings
of $g_i$ (see figure 4).

Let $K$ be the union of all those graphs : $K = \bigcup^n_{i=1}
K_{T_i}$. Some of the components of $K$ have edges with one vertex in a hole of a
 slice decomposition. Let $K'$ be the graph obtained from $K$ when one has
removed all these components.

There are two kind of connected components of $K'$ : the components $K_i$ for
which a small tubular neighborhood $NK_i$ is such  that $NK_i\setminus K_i$ is
disconnected (type I), and those for which it is connected (type II).

\subsection{$G$ as a graph of groups.}

We now split $G$ in a graph of groups by cutting $P$ along the graph $K$.

The graph of groups we consider is as follows. Its vertices are of two
kinds. First there are the connected
components of $P\setminus  K'$, and the groups are the fundamental groups of
those components.
There are also the components of $K'$ of type II, and the groups are the fundamental
groups of a small tubular neighborhood. The  edges of  the graph of groups are
the components $K_i$ of $K'$, and their groups are either $\pi_1(NK_i)$, the
fundamental group of a small tubular neighbourhood, in the case of a component
of type I, or $\pi_1(NK_i \setminus K_i)$ otherwise, in type II. Note that in
this case, $\pi_1(NK_i \setminus K_i)$ is of index two in $\pi_1(NK_i)$.

\begin{lemma}(\cite{Del} Lemma III.2.b)

Let $H$ be a subgroup of $G$ stabilizing an edge of the graph of groups. Then
  $h(H)$ is a subgroup of $\Gamma$ that has an orbit in $\mathcal{K}$ which is
contained in a slice. In particular, this orbit is finite. \end{lemma}

For the proof, see \cite{Del}.

 In the case of hyperbolic groups, one deduces that the subgroup is finite ; in our case,
 by Lemma \ref{lem;virtparab}, it is either finite or parabolic.

\begin{coro}\label{coro;leaflambda}
  If the map $h$ has no accidental parabolic, then the graph of
groups is a trivial splitting, and  $h(G)$ is the image of a vertex group
corresponding to a leaf $\lambda$ in $P$, containing singular points of the
lamination : $h(G)$ is conjugated to the image of $\pi_1(\lambda)$
(only defined up to conjugacy).
\end{coro}

\subsection{If $h$ has no accidental parabolic.}

In all the following, we assume that $h$ has no accidental parabolic: we
can apply Corollary \ref{coro;leaflambda}.

Let $P_\Gamma$ be a Van Kampen polyhedron for $\Gamma$, for a presentation with a finite generating
set : it is a cell complex of dimension 2, whose 1-skeleton
consists of finitely many loops.
The set of vertices of its
universal cover is $\widetilde{P_\Gamma}^0$, and after the choice of a base point, we identify it with $\Gamma$, 
hence also with $\Gamma p =\{\gamma p , \gamma \in \Gamma \} \subset \mathcal{K}$, the set of vertices of finite valence of
 $\mathcal{K}$:  $\widetilde{P_\Gamma}^0 \approx \Gamma \approx (\Gamma p)$.

\subsubsection{Lifting slices of $\mathcal{K}$ in $\Gamma$}

From the identification above, we have an inclusion of $\widetilde{P_\Gamma}^0$ in $\mathcal{K}$.  
 We want to define good pre-images 
in $\widetilde{P_\Gamma}^0 \approx \Gamma$ of
slices of cylinders in $\mathcal{K}$. We define the pre-image $S_\Gamma$ of an arbitrary slice $S$ as follows.

 If $S$ is a regular slice of a cylinder in $\mathcal{K}$, which is \emph{not}
reduced to a vertex of infinite valence, then we say that $S_\Gamma$ is $S_\Gamma = \{ \gamma \in \Gamma \, | \;
\exists s_1, s_2 \in S, \, |s_1 - \gamma p | + |\gamma p -s_2|=|s_1 - s_2 |\}$, the set
of elements of $\Gamma$ that send the base point $p$ of $\mathcal{K}$ on a
geodesic segment with ends in $S$.

If $S=\{v\}$ is a parabolic slice of a cylinder in $\mathcal{K}$,
 or a regular slice reduced to a point of infinite valence,
 then we define $S_\Gamma$ to be the set $\{\gamma \in \Gamma \, | \, |\gamma p -v |
=1 \}$. It is a coset of a parabolic subgroup of
$\Gamma$.

\subsubsection{The map  $\tilde{h}: \tilde{P} \to \widetilde{P_\Gamma}$ }

Let $\tilde{P}$ be the universal cover of the polyhedron $P$, and $*$ a base point in it.
For every $i = 1\dots k$, for every edge $c_i$ in the one skeleton of $P$, we
denote by $\tilde{c}_i$ its image in $\tilde{P}$ starting at $*$. We define markings on $\tilde{c}_i$ such that they map exactly on the ones of $c_i$ by the quotient map, and we extend the construction by $G$-equivariance in $\tilde{P}$. Every edge of the
1-skeleton of $\tilde{P}$ is hence marked by consecutive markings.

Recall that $\widetilde{P_\Gamma}$ is the universal cover of $P_\Gamma$. 
 The
morphism $h$ can be realized as a continuous $G$-equivariant map 
$\tilde{h}:\tilde{P} \to \widetilde{P_\Gamma}$
 such that for all $i=1\dots k$,
$\tilde{h}(\tilde{c}_i)$ is a path from 
$\tilde{h}(*)$ to $h(g_i)\tilde{h}(*)$,
where $g_i$ denotes the element of $G$ associated to $c_i$.

The map $\tilde{h}$ is completely defined on the vertices of $\tilde{P}$. We now choose the images of the markings of each edge $\tilde{c_i}$ ($i=1\dots k$). 

First, if $m_i^j$ is any marking of $c_i$ associated to a slice $S$ 
(without restriction), 
 and if $\tilde{m}_i^j$ is its image in $\tilde{c}_i$ then
$\tilde{h}(\tilde{m}_i^j)$ is equal to a vertex $\gamma \tilde{h}(*)$ of
$\widetilde{P_\Gamma}$, such that $\gamma \in S_\Gamma$.

Second, if $m_i^j$ is a marking of $c_i$ associated to a parabolic slice $S$,
 there is an unique marking adjacent to $m_i^j$ in $c_i$, which is
associated to a slice $S' \neq S$. Then we require that
$\tilde{h}(\tilde{m}_i^j)=\gamma \tilde{h}(*)$, where $\gamma \in S_\Gamma$ is
such that $\gamma p$ lies on some geodesic from $v$ to a point of
$S'$ in $Cyl(p,h(g_i)p)$. 
We denote by $S_\Gamma (i,j)$ the set of such elements 
$\gamma \in S_\Gamma$.
 Note that the images of the two markings of a parabolic slice
might be very far from each other in $\Gamma$, in the same coset of parabolic
subgroup.

Third, if $m_i^j$ is a marking of $c_i$ associated to a regular slice $S$
reduced to a vertex of infinite valence $S=\{v\}$, then we require that
$\tilde{h}(\tilde{m}_i^j)=\gamma \tilde{h}(*)$, where $\gamma \in S_\Gamma$ is
such that $\gamma p$ lies on some geodesic from $v$ to a point of a slice
adjacent to $S$ in $Cyl(p,h(g_i)p)$. We denote by $S_\Gamma (i,j)$ the set of such elements
$\gamma \in S_\Gamma$.

We can assume that $\tilde{h}(\tilde{c}_i)$ is a geodesic between the images
of consecutive markings, but this is not essential.

\begin{lemma}\label{lem;aidepourvertex}
 Let $v$ ba a vertex of $\mathcal{K}$ of infinite valence, such that $\{v\}$ is a slice (either parabolic or regular) of the cylinder $Cyl(p,h(g_i)p)$.  Let $m_i^j$ be the  marking on $c_i$, associated to the slice $S=\{v\}$. 
The diameter of $S_\Gamma (i,j)$ in $\Gamma$ (for the word metric) 
is at most $2000\delta(2\Theta+1)$.

 \end{lemma}

{\it Proof : }  Let $\gamma_1$ and $\gamma_2$ be in $S_\Gamma (i,j)$.
 There are points $v_1$
and $v_2$ in  slices $S'_1$ and $S'_2$ adjacent to $S$ in $Cyl(p,h(g_i)p)$.
 By Corollary \ref{coro;slicesareclose}, $|v-v_i|\leq  1000\delta$, and for some
geodesic segments, $\MaxAng([v,v_i])\leq 2\Theta$, for $i=1,2$.

First assume that $S$ is a parabolic slice. Then, by the definition of $S_\Gamma
(i,j)$, $S'_1 = S'_2$. By Lemma \ref{lem;pointes}, $\Ang_v([v,v_1],[v,v2]) \leq
14D\leq \Theta $. Therefore, by the majoration (\ref{eq;majoration}), 
we can deduce that $|\gamma_1^{-1} \gamma_2| \leq 2000\delta(2\Theta+1)$.

Secondly assume that $S$ is a regular slice. Then there is no parabolic slice
between $S'_1$ and $S'_2$. By Lemma \ref{lem;noanglewithoutparabolic}, $\Ang_v([v,v_1],[v,v2]) \leq
2\Theta$. Therefore, again by the majoration (\ref{eq;majoration}), $|\gamma_1^{-1} \gamma_2| \leq  2000\delta(2\Theta+1) $. $\square$

\subsubsection{Bounding the lengths of the images of leaves of $\Lambda$ in
$P_\Gamma$}

The equivariant map $\tilde{h}$ induces a continuous map $h:P\to P_\Gamma$.

The next lemma is an analogue of Lemma II.1 in \cite{Del}, but cannot be
deduced from it, because of the presence of
parabolic slices.

\begin{lemma}\label{lem;bound_reg_leaf}

   Let  $l_1,\dots,l_m$ be  a sequence of regular arcs of $\Lambda$,  where
$l_i$ links the marking $\iota(l_i)$ to the marking $\tau(l_i)$, and where
$\tau(l_i)=\iota(l_{i+1})$.
If the path $l_1l_2...l_m$ has no loop, then the
path $h(l_1l_2...l_m)$ in $P_\Gamma$ is homotopic, with fixed ends, 
to a path in
the 1-skeleton of $P_\Gamma$, of length less than $20 000\delta(\Theta
+1)\times n$ (for the graph metric of the 1-skeleton).

\end{lemma}

{\it Proof : } As the arcs are all regular, all the markings 
involved are associated to the
same slice of $\mathcal{K}$, say $S$. Let us lift the path $l_1l_2...l_k$ in a
path $\widetilde{l_1l_2...l_k}$ of $\tilde{P}$, starting at the marking
$\tilde{m}_i^j$,  where   $m_i^j= \iota(l_1)$.  Thus, this path is mapped in 
$\widetilde{P_\Gamma}$ on a path that stays in $S_\Gamma$. As
$\widetilde{P_\Gamma}$ is simply connected, this path is homotopic to any path
in the 1-skeleton that has the same ends.

There are two main cases to study, namely if the slice is regular not reduced
to a single point of infinite valence, or if it is reduced to a single point of
infinite valence (including the case of parabolic slices). If the second
case, we will have to discuss whether an
adjacent arc of the lamination is regular or not.

 First, if the slice $S$ is regular, not reduced to a parabolic point, then the
end points $v_0$ and $v_m$ of $\tilde{h}(\widetilde{l_1l_2...l_m})$ 
are vertices
of the form $v_0=\gamma_0 \tilde{h}(*)$ for $\gamma_0 \in S_\Gamma$, and
$v_m=\gamma_m \tilde{h}(*)$ for $\gamma_m \in S_\Gamma$. Therefore, there exist
$s_0$ and $s'_0$ in $S$ and a geodesic segment  $[s_0,s'_0]$ in $\mathcal{K}$
containing $\gamma_0 p$ (and similarly for $\gamma_m$).   By Lemma
\ref{lem;slicesaresmall}, we have a path  from  $\gamma_0 p$ to  $\gamma_m
p$ of length at most $3\times 200\delta$, and of maximal angle at most
$2\Theta$. Therefore, by the majoration (\ref{eq;majoration}), the distance 
in the 1-skeleton of $\widetilde{P_\Gamma}$
between $v_0$ and $v_m$ is at most $600\delta(2\Theta +1)$.

Secondly, we assume that $S$ is a parabolic slice
or a regular slice reduced to a single vertex of infinite valence.
 Then in the edge containing the marking
$\iota(l_i)$, there is one (and only one, if the slice is parabolic) marking
$m_{\iota,i}$ adjacent to $\iota(l_i)$ that is not associated to $S$. In the
edge containing the marking $\tau(l_i)$, there is only one marking $m_{\tau,i}$
adjacent to $\tau(l_i)$ that is not associated to $S$, and that is linked to
$m_{\iota,i}$ by an arc (regular or singular) of the lamination of the triangle
or digon. These markings are associated to regular slices (cf Remark 3).

There are two possibilities.

 In the triangle containing $l_i$, it is possible that
$[m_{\iota,i}, m_{\tau,i}]$ is a regular arc of $\Lambda$.
 Let $l_{i_0} \dots l_{i_q}$ a maximal subpath such that this property holds
at each step.  By Lemma
\ref{lem;aidepourvertex},  the end points of the image of $l_{i_0} \dots
l_{i_q}$ in $\widetilde{P_\Gamma}$ are at distance at most
$2000\delta(2\Theta+1)$ in the
1-skeleton of $\widetilde{P_\Gamma}$. Therefore, the image of $l_{i_0} \dots
l_{i_q}$ in $\widetilde{P_\Gamma}$ is homotopic  with fixed ends,
to a path in the 1-skeleton of length less than $2000\delta(2\Theta+1)$.

Assume now that $[m_{\iota,i}, m_{\tau,i}]$ is not a regular arc of $\Lambda$.
 That is that   $l_i$    is  one of the three  regular leaf of a triangle
that is adjacent to a singular leaf. Note that in a path $l_1\dots l_m $ without
loop, this can only happen $3n$ times, where $n$ is the number of triangles.

Let $S$ be the slice of the cylinders of the triangle containing
$l_i$, associated to $\iota(l_i)$ and $\tau(l_i)$.
Let $S_\iota$ be the slice associated to $m_{\iota,i}$, and $S_\tau$ be the
slice associated to $m_{\tau,i}$.

  In order to bound the distance between the images of $\iota(l_i)$
and $\tau(l_i)$, it is enough to bound the maximal angle of geodesics between
elements of $S_\iota$ and $S_{\tau}$. Let $v_{\iota }$ be in $S_\iota$,
$v_{\tau }$ be in $S_{\tau}$.

 We claim that, given a geodesic segment
$[v_{\iota }, v_{\tau }]$ in $\mathcal{K}$, its maximal angle is at most
$5\Theta$.

If $S$ is a regular slice, it is the triangular inequality for angles in the two
edges of the triangle sharing $S$.

  If $S$ is parabolic, we consider a segment between
 $v_{\iota }$ and $v_{\tau }$ that passes through the vertex of the
slice $S$. By Lemma \ref{lem;slicesaresmall}, it has no angle larger than
$2\Theta$ except possibly at $S$, and if its angle is larger than $5\Theta$ at
this point, $S$ would be a parabolic slice of the third side of the triangle. By
the construction of the leaves in a triangle, the marking $\tau(l_i)$ should be
on this side, which is not the case.

  Therefore, the distance between the images of $\iota(l_i)$ and $\tau(l_i)$ in
the 1-skeleton of $\widetilde{P_\Gamma}$ is at most $5\Theta
$.

   For a path $l_1l_2...l_m$ without
loop, such a situation can happen only $3n$ times, where $n$ is the number of
triangles. Therefore the distance between the endpoints of its image,  in
the 1-skeleton of $\widetilde{P_\Gamma}$, is at most
$3n\times(2000\delta(2\Theta+1) + 5 \Theta) + 2000\delta(2\Theta+1)$. This is
less than  $20 000\delta(\Theta+1)\times n$. $\square$

\vskip .3cm

\begin{lemma}\label{lem;bound_sing_leaf}

An arc of $\Lambda$ linking two markings corresponding to slices in a hole of a
same triangle, maps on a path which is homotopic, with fixed ends, to a path in
the 1-skeleton of $P_\Gamma$, of length less
than $(\varphi(n)+1)\times (40 000\delta(\Theta+1))$.
\end{lemma}

{\it Proof : }  Such  an arc  is homotopic with
 fixed ends in $P$ to a path tracking back on the first side of the triangle,
 until the first regular arc to the other side, and then tracking on this side
 to the suitable marking. By theorem \ref{theo;slices}, this path enters in at
most  $2\times (10\varphi(n)+1)$ slices, none of them having an angle
superior to $5\Theta$. Therefore, by the majoration (\ref{eq;majoration}), 
the distance between the end points of the image is inferior to  
$2\times (10\varphi(n)+1)\times (1000\delta(2\Theta +1))$ in the
1-skeleton of the universal cover of $P_\Gamma$. $\square$

\subsubsection{Image of the leaf $\lambda$}

We need a lemma from \cite{Del}.

\begin{lemma}(\cite{Del} Lemma III.4)

 Let $L$ be a connected graph, $L_1$ be its 1-skeleton, and $E$ a  metric
space. Let $h:L \to E$ be a continuous map. Let $E'$ be a subset of $E$.
Assume that :

  1) For all edge $l$ in $L_1$, $h(l)$ is homotopic in $E$, with fixed ends, to
a curve in $E'$ of length less than the constant $M$.

  2) There exists a finite set of edges $L'_1 \subset L_1$ such that a path
without loop,  made of consecutive edges $l_1, \dots, l_k$ in $L_1\setminus
L'_1$, has its image by $h$ homotopic in $E$ (with fixed ends) to a curve
 in $E'$ of length less than $M$.

 Then, for all vertex $s$ of $L$, $h_* (\pi_1(\Lambda,s))$ is generated by
curves in $E'$ of length inferior to $(4\Card(L'_1) +3)\times  M$. \end{lemma}

{\it Proof : } Let $T$ be a maximal tree in $L$.
The group $h_* (\pi_1(\Lambda,s))$ is generated by the images of the loops of
the form $[s,s']e[s",s]$, where the segments $[s,s'] $ and $[s",s]$ are in
$T$, and where $e$ is an edge from $s'$ to $s"$ in  $L$. In particular, the paths  $[s,s'] $ and
$[s",s]$ do not contain any loop, and contain at most $\Card(L'_1)$ edges of
$L'_1$. Each of those two segments are the concatenation of at most
$\Card(L'_1)+1$ segment without loop made of consecutive edges in
$L_1\setminus L'_1$, with at most  $\Card(L'_1)$ edges of
$L'_1$. Therefore the image of  $[s,s']$ by $h$ is homotopic in $E$, with fixed
ends, to a curve in $E'$ of length less than $(2\Card(L'_1) +1)\times  M$, and
the same is true for the image of $[s",s]$. Finally, the image of the edge $e$
is homotopic with fixed ends to a curve of $E'$ of length at most $M$, this
gives the result. $\square$

Finally, we can prove Theorem \ref{theo;acci}. Given a morphism $h :G\to
\Gamma$ without accidental parabolics, we set $E=P_\Gamma$,  $E'$ its
1-skeleton, and $L=\lambda$, the singular leaf of $\Lambda$ given by Corollary
\ref{coro;leaflambda}. We choose  $L'_1$ to be the set of arcs joining two
markings of a hole of a triangle, via the singular point of this triangle, and
$M =40 000\delta(\varphi(n)+1)(\Theta+1) \}$  (which is superior to $
20 000\delta (\Theta +1) \times n$). By Lemma \ref{lem;bound_reg_leaf} and Lemma
\ref{lem;bound_sing_leaf}, the assumptions of the previous lemma are fulfilled.
We get that $h(G)$ is conjugated to a subgroup of $\Gamma$ generated by curves
in the 1-skeleton of $P_\Gamma$ of length bounded by $(4\times n\times
(30\varphi(n))^2 +3)\times M$. There are finitely many such curves. Hence, there
are finitely many such subgroups, therefore this implies Theorem
\ref{theo;acci}. $\square$

\section{Appendix : Coarse-piecewise-geodesics are $\lambda$-quasi-geodesics.}

In this appendix, we give a simple proof that coarse piecewise geodesics (Definition \ref{def;cpg}) are $\lambda$-quasi-geodesics (Proposition \ref{prop;final}). Let $\mathcal{K}$ be an hyperbolic graph, and $l$ a constant greater than $\mu$ (see section 2 for the constants). Let $f:[a,b]\to \mathcal{K}$ be a coarse piecewise geodesic, for the subdivision of $[a,b]$ : $a=c_1\leq d_1\leq \dots \leq c_n \leq d_n =b$.

\begin{lemma} \label{lem;antefinal}

Let $i$ be an integer in $[1,n]$. Let $t\in [c_i,d_i]\subset [a,b]$ be such that  $|t-c_i|\geq 4\epsilon $ and $|t-d_i|\geq 4\epsilon$. Then $f(t)$ is at distance at most $2\delta$ from a geodesic segment $[f(a), f(b)]$. 

\end{lemma}

Let us mention that the proof is similar to the one of Lemma \ref{lem;2delta}.

{\it Proof : } As $f|_{[c_i,d_i]}$ is  a $\mu$-local geodesic, the restriction  $f|_{[(t-4\epsilon), (t+4\epsilon)]}$ is a geodesic segment whose ends are at distance at most $2\epsilon$ from a segment $[f(a),f(b)]$. Let $w_1$ and $w_2$ be points in this segment realizing the minimal distance to $f(t-4\epsilon)$, $f(t+4\epsilon)$. The triangle $(f(t-4\epsilon),w_1,f(t+4\epsilon) )$ and $(w_1,w_2,f(t+4\epsilon) )$ are $\delta$-thin, therefore, $v$ is at distance at most $2\delta$ from $[w_1,w_2]$. $\square$

\begin{prop} \label{prop;final}
 
 Let $t_1$ and $t_2$ be such that $a\leq t_1< t_2 \leq b$. 
Then $|f(t_1)-f(t_2)| \geq \frac{1}{\lambda} |t_1-t_2|$.

\end{prop}

{\it Proof : } Either there is a number $u_1$ such that $|u_1 - t_1|\leq 5 \epsilon$ and such that $u_1$ satisfies the hypothesis of Lemma \ref{lem;antefinal}, or $|a - t_1|\leq 5 \epsilon$ (in this case we write $u_1 =a$). In both cases, $f(u_1)$ is at distance at most $2\delta$ from a point $v_1$ in $[f(a),f(b)]$. 
Let $k$ be a positive integer such that $t_1 + 1000 k\lambda\delta \leq t_2$. 
Then, there exists $u_{k+1}$ a number such that $|t_1 + 1000 k\lambda\delta - u_{k+1}| \leq 5 \epsilon$, and satisfying the hypothesis of Lemma \ref{lem;antefinal}. Therefore there exists $v_{k+1}$ on $[f(a),f(b)]$ at distance at most $2\delta$ from $f(u_{k+1})$.
 
Let $m$  be the maximal number such that $t_1 + 1000 m\lambda \delta \leq t_2$.

By definition of coarse-piecewise-geodesics, for all $k \in [1,m+1]$, $f|_{[u_k,u_{k+1}]}$ is a $\frac{\lambda}{2}$-quasi-geodesic. Therefore, $|f(u_k)-f(u_{k+1})| \geq  \frac{2}{\lambda}|u_{k+1} - u_k|$. We deduce that $|v_k-v_{k+1}|\geq   \frac{2}{\lambda}|u_{k+1} - u_k| - 4\delta$. Therefore, by summing, $|v_1 - v_{m+1}| \geq   \frac{2}{\lambda} \times |u_{m+1} - u_1| - 4m\delta $. 

Moreover,  $|v_{m+1} - f(t_2)|  \leq 5\epsilon + 1000\lambda\delta +2\delta$, and $|v_{1} - f(t_1)|  \leq 5\epsilon +2\delta $. Therefore, $|f(t_1)-f(t_2)| \geq  \frac{2}{\lambda} \times |u_{m+1} - u_1| - 4(m+1)\delta - 10\epsilon - 1000\lambda\delta$. Since $|u_0-t_1|\leq  5\epsilon$ and $|u_{m+1}-t_2|\leq 5\epsilon + 1000\lambda\delta$, we get that $|f(t_1)-f(t_2)| \geq \frac{2}{\lambda} \times (|t_2 - t_1| - 10 \epsilon - 1000 \lambda\delta)  - 4(m+1)\delta - 10\epsilon - 1000\lambda\delta$. 
Since $m \leq \frac{|t_2 - t_1|}{1000 \lambda\delta}$ and since $\frac{\lambda}{2} \geq 1$,  we deduce that  
$|f(t_1)-f(t_2)| \geq (\frac{2}{\lambda}-\frac{4\delta}{1000\lambda\delta})\times |t_2 - t_1| - 20\times  (\epsilon + 100 \lambda \delta)$.

Finally one has $|f(t_1)-f(t_2)| \geq \frac{1}{\lambda}|t_2 - t_1| \times (2-\frac{1}{250}) + 20 ( \epsilon + 100 \lambda\delta)$. If 
$|t_2-t_1| \geq 40 \lambda \times (\epsilon + 100 \lambda\delta)$, then $|f(t_1)-f(t_2)| \geq \frac{1}{\lambda}|t_2 - t_1|$. Otherwise, the result comes from the assumption that $f$ is a  local  quasi-geodesic. $\square$

\end{document}